\documentclass{amsart}

\usepackage{amssymb,amsxtra}

\numberwithin{equation}{section}
\newcommand\note[1]%
{$^\dagger$\marginpar{\sf\footnotesize{$^\dagger${#1}}}}

\def\today{\number\day\space\ifcase\month\or January\or February\or
March\or April\or May\or June\or July\or August\or September\or
October\or November\or December\fi\space\number\year}

\swapnumbers
\newtheorem{theorem}[equation]{Theorem}
\newtheorem{proposition}[equation]{Proposition}
\newtheorem{lemma}[equation]{Lemma}
\newtheorem{corollary}[equation]{Corollary}

\theoremstyle{definition}
\newtheorem{definition}[equation]{Definition}
\newtheorem{example}[equation]{Example}
\newtheorem{remark}[equation]{Remark}

\newcommand\eu{\mathfrak}
\newcommand\lie{\mathfrak}
\renewcommand\k{\lie{k}} 
\renewcommand\t{\lie{t}}
\newcommand\g{\lie{g}}

\newcommand\z{\lie{z}}

\newcommand\bb{\mathbb}

\newcommand\Z{\bb{Z}} 

\newcommand\R{\bb{R}} 
\newcommand\C{\bb{C}}
\renewcommand\H{\bb{H}}
\renewcommand\P{\bb{P}}

\newcommand\ca{\mathcal}

\newcommand\Ad{\operatorname{Ad}}
\newcommand\ad{\operatorname{ad}}

\renewcommand\star{\operatorname{star}}

\newcommand\RR{\operatorname{RR}}

\newcommand\Ind{\operatorname{Ind}}

\newcommand\pr{\operatorname{pr}}

\newcommand\Hom{\operatorname{Hom}}
\newcommand\Spec{\operatorname{Spec}}

\newcommand\lk{\operatorname{lk}}
\newcommand\slk{\operatorname{slk}}

\newcommand\id{\operatorname{id}}

\renewcommand\Im{\operatorname{Im}}

\newcommand\Lie{\operatorname{Lie}} 

\newcommand\SO{\operatorname{\mathbf{SO}}}

\newcommand\SL{\operatorname{\mathbf{SL}}}
\newcommand\SU{\operatorname{\mathbf{SU}}}
\newcommand\U{\operatorname{\mathbf{U}}}

\newcommand\abs[1]{\lvert#1\rvert}

\newcommand\norm[1]{\lVert#1\rVert}

\newcommand\quot[1][\kern.3ex]{/\kern-.7ex/_{\kern-.4ex#1}}
\newcommand\bigquot[1][\,\,]{\big/\kern-.85ex\big/_{\!\!#1}}

\newcommand\powl{[\kern-.3ex[}
\newcommand\powr{]\kern-.3ex]}
\newcommand\bigpowl{\bigl[\kern-.6ex\bigl[}
\newcommand\bigpowr{\bigr]\kern-.6ex\bigr]}

\newcommand\longhookrightarrow{\lhook\joinrel\longrightarrow}

\newcommand\sur{\longrightarrow\kern-1.9ex\to}
\newcommand\iso{\longhookrightarrow\kern-1.9ex\to}

\newcommand\cross{\Phi\inv(\t^*_+)}

\newcommand\antiddots{\mathinner{%
\mkern1mu\raise1pt\vbox{\kern7pt\hbox{.}}
\mkern2mu\raise4pt\hbox{.}\mkern2mu\raise7pt\hbox{.}\mkern1mu}}

\newcommand\inv{^{-1}} 

\renewcommand\subset{\subseteq}

\newcommand\st{^{\mathrm{s}}}
\newcommand\sst{^{\mathrm{ss}}}

\newcommand\reg{_{\mathrm{reg}}}
\newcommand\prin{_{\mathrm{prin}}}

\newcommand\impl{_{\mathrm{impl}}}

\begin{document} 


\title{Symplectic Implosion}

\author{Victor Guillemin}

\address{Department of Mathematics, Massachusetts Institute of
Technology, Cambridge, Massachusetts 02139-4307}

\email{vwg@math.mit.edu}

\author{Lisa Jeffrey}

\address{Department of Mathematics, University of Toronto, Toronto,
Ontario M5S 3G3}

\email{jeffrey@math.toronto.ca}

\author{Reyer Sjamaar}

\address{Department of Mathematics, Cornell University, Ithaca, New
York 14853-7901} 

\email{sjamaar@math.cornell.edu}

\thanks{V. Guillemin was partially supported by NSF Grant DMS-9625714,
L. Jeffrey by an Alfred P. Sloan Research Fellowship and by an NSERC
Grant, and R. Sjamaar by an Alfred P. Sloan Research Fellowship and by
NSF Grant DMS-0071625}

\date{18 January 2001}


\begin{abstract}
Let $K$ be a compact Lie group.  We introduce the process of
symplectic implosion, which associates to every Hamiltonian
$K$-manifold a stratified space called the imploded cross-section.  It
bears a resemblance to symplectic reduction, but instead of
quotienting by the entire group, it cuts the symmetries down to a
maximal torus of $K$.  We examine the nature of the singularities and
describe in detail the imploded cross-section of the cotangent bundle
of $K$, which turns out to be identical to an affine variety studied
by Gelfand, Vinberg, Popov, and others.  Finally we show that
``quantization commutes with implosion''.
\end{abstract}

\maketitle

\tableofcontents


\section{Introduction}

According to Cartan and Weyl, a finite-dimensional representation of a
compact connected Lie group $K$ is determined up to isomorphism by its
highest-weight vectors.  In the language of the orbit method, the
classical analogue of a representation is a symplectic manifold $M$
equipped with a Hamiltonian action of the group $K$.  The classical
analogue of the collection of highest weights is then the moment, or
Kirwan, polytope of $M$.  This paper deals with the question, what is
the classical analogue of the set of highest-weight vectors?

In answer to this question we construct a space called the
\emph{imploded cross-section} of $M$.  It is defined by fixing a Weyl
chamber of $K$, taking its inverse image under the moment map, and
quotienting the resulting subset of $M$ by a certain equivalence
relation.  While this subquotient construction is reminiscent of
symplectic reduction,
the imploded cross-section is almost always singular, whereas
symplectic quotients often are not.  For example, the imploded
cross-section of the cotangent bundle $T^*K$ is singular unless the
commutator subgroup of $K$ is a product of copies of $\SU(2)$.  It is
not even an orbifold unless the universal cover of $[K,K]$ is a
product of copies of $\SU(2)$.  (See Section \ref{section;kaehler}.)
The singularities are however not completely arbitrary.  Like singular
symplectic quotients, imploded cross-sections stratify naturally into
symplectic manifolds and the structure of the singularities is locally
conical.  The imploded cross-section is defined in Section
\ref{section;construction} and its stratification is studied in
Section \ref{section;stratification}.

Moreover, there is a residual action on $M\impl$ of a maximal torus
$T$ of $K$, which is Hamiltonian in a suitable sense.  The image of
the $T$-moment map on $M\impl$ is equal to the moment polytope of $M$.
The classical analogue of the Cartan-Weyl theorem is then the
following assertion: for each $\lambda$ in the Weyl chamber the
quotients $M\quot[\lambda]K$ and $M\impl\quot[\lambda]T$ are
symplectomorphic.  (The notation $\quot[\lambda]$ stands for
symplectic reduction at level $\lambda$.)  In this way the process of
symplectic implosion abelianizes Hamiltonian $K$-manifolds at the cost
of introducing singularities.  This issue is dealt with in Section
\ref{section;abelian}.

The imploded cross-section of $T^*K$ enjoys two special properties.
The first, which is explored in Section \ref{section;universal}, is
that $(T^*K)\impl$ carries an additional $K$-action which commutes
with the $T$-action, and that the imploded cross-section of any
Hamiltonian $K$-manifold $M$ can be obtained as a symplectic quotient
of the product $M\times(T^*K)\impl$ with respect to $K$.  For this
reason we call $(T^*K)\impl$ the \emph{universal} imploded
cross-section.  The second property, which is investigated in Section
\ref{section;kaehler}, says that $(T^*K)\impl$ possesses the structure
of an irreducible complex affine variety and that the symplectic
structure is K\"ahler.  The $K$-action extends to an algebraic action
of the complexified group $G=K^\C$ and the $G$-orbits are identical to
the symplectic strata.  The open orbit is of type $G/N$, where $N$ is
a maximal unipotent subgroup of $G$.  In fact, $(T^*K)\impl$ can be
identified with the basic affine variety introduced by Bernstein et
al.\ \cite{bernstein-gelfand-gelfand;differential-operators-base}.
Thus implosion is the symplectic counterpart of taking the quotient of
a $G$-variety by the action of $N$.

The result of Section \ref{section;quantization} reinforces further
the analogy between imploded cross-sections and highest-weight
vectors.  It asserts that ``quantization commutes with implosion'' in
the following sense.  Assuming that $M$ is prequantizable, we define
its quantization to be the $K$-equivariant Riemann-Roch number
$\RR(M)$ with coefficients in the prequantum line bundle.  The
quantization of $M\impl$ is defined as the $T$-equivariant
Riemann-Roch number of a suitable desingularization.  The result is
that $\RR(M\impl)$ is equal as a virtual $T$-module to the space of
highest-weight vectors in $\RR(M)$.  


Many of the results in Sections
\ref{section;construction}--\ref{section;kaehler} are taken from an
unpublished manuscript dating from 1993.  They have recently found an
application in the theory of vector bundles on Riemann surfaces (see
Hurtubise and Jeffrey
\cite{hurtubise-jeffrey;representations-weighted}), which is why we
are making available this updated and expanded version.  More
precisely, this application involves symplectic implosion for
group-valued moment maps, which we hope to take up in a sequel to this
paper.

\section{The construction}\label{section;construction}

Let $(M,\omega)$ be a connected symplectic manifold and let $K$ be a
compact connected Lie group acting on $M$ in a Hamiltonian fashion
with equivariant moment map $\Phi\colon M\to\k^*$, where $\k=\Lie K$.
Our sign convention for the moment map is
$d\Phi^\xi=\iota(\xi_M)\omega$, where $\xi_M$ denotes the vector field
on $M$ induced by $\xi\in\k$, and $\Phi^\xi=\langle\Phi,\xi\rangle$
the component of the moment map along $\xi$.

We choose once and for all a maximal torus $T$ of $K$ and a closed
fundamental Weyl chamber $\t^*_+$ in the dual of the Cartan subalgebra
$\t=\Lie T$.  The Weyl chamber is the disjoint union of $2^r$ open
faces (sometimes called walls), where $r$ is the rank of the
commutator subgroup $[K,K]$.  The \emph{principal face} $\sigma\prin$
for $M$ is the minimal face $\sigma$ such that
$\Phi(M)\cap\t^*_+\subset\bar\sigma$.  In most cases of interest it is
equal to $(\t^*_+)^\circ$, the interior of the Weyl chamber.  The
symplectic cross-section theorem (see below) says that
$\Phi\inv(\sigma\prin)$ is a $T$-stable symplectic submanifold of $M$.
We want to ``complete'' this submanifold by adding lower-dimensional
symplectic strata.  An obvious guess is to take
$\cross=\Phi\inv(\bar\sigma\prin)$, but to turn this into a symplectic
object, we need to contract it along certain directions in the
boundary components in the following manner.  Two points $m_1$ and
$m_2$ in $\cross$ are \emph{equivalent}, $m_1\sim m_2$, if there
exists $k\in[K_{\Phi(m_1)},K_{\Phi(m_1)}]$ such that $m_2=km_1$.  By
equivariance of the moment map, $m_1\sim m_2$ implies
$\Phi(m_1)=\Phi(m_2)$, so $\sim$ is an equivalence relation.

\begin{definition}\label{definition;implosion}
The \emph{imploded cross-section} of $M$ is the quotient space
$M\impl=\cross/{\sim}$, equipped with the quotient topology.  The
quotient map $\cross\to M\impl$ is denoted by $\pi$.  The
\emph{imploded moment map} $\Phi\impl$ is the continuous map
$M\impl\to\t^*_+$ induced by $\Phi$.
\end{definition}

The image of $\Phi\impl$ is equal to $\Phi(M)\cap\t^*_+$.  All points
in a face $\sigma$ have the same centralizer, denoted $K_\sigma$, and
therefore $M\impl$ is set-theoretically a disjoint union of orbit
spaces,
\begin{equation}\label{equation;decomposition}
M\impl=\coprod_{\sigma\in\Sigma}\Phi\inv(\sigma)/[K_\sigma,K_\sigma].
\end{equation}
Here $\Sigma$ denotes the collection of faces of $\t^*_+$.  We define
a partial order on $\Sigma$ by putting $\sigma\le\tau$ if
$\sigma\subset\bar\tau$.

\begin{lemma}\label{lemma;hausdorff}
The projection $\pi$ is proper and $M\impl$ is Hausdorff\upn, locally
compact\upn, and second countable.  If $M$ is compact\upn, then so is
$M\impl$.  Each subspace in the decomposition
\eqref{equation;decomposition} is locally closed in $M\impl$.
\end{lemma}

\begin{proof}
First we show that $\pi$ is closed.  Let $C\subset\cross$ be closed.
We need to show that
$$
\pi\inv(\pi(C))=\coprod_{\sigma\in\Sigma}\Phi\inv(\sigma)
\cap[K_\sigma,K_\sigma]\cdot C
$$
is closed.  Let $\{m_i\}$ be a sequence in $\pi\inv(\pi(C))$
converging to $m\in\cross$.  Let $\sigma$ be the face containing
$\Phi(m)$.  After passing to a subsequence we may assume that all
$\Phi(m_i)$ are in the same face, say $\tau$.  Then $\sigma\le\tau$,
so $K_\tau\subset K_\sigma$ and
$[K_\tau,K_\tau]\subset[K_\sigma,K_\sigma]$.  Since $[K_\tau,K_\tau]$
is compact, $[K_\tau,K_\tau]\cdot C$ is closed, so
$m\in[K_\tau,K_\tau]\cdot C$.  Hence
$$
m\in\Phi\inv(\sigma)\cap[K_\tau,K_\tau]\cdot
C\subset\Phi\inv(\sigma)\cap[K_\sigma,K_\sigma]\cdot
C\subset\pi\inv(\pi(C)),
$$
i.e.\ $\pi\inv(\pi(C))$ is closed.  The fact that $\pi$ has compact
fibres now implies that it is proper.  The stated properties of
$M\impl$ follow easily, and the local closedness of
$\Phi\inv(\sigma)/[K_\sigma,K_\sigma]=\pi\bigl(\Phi\inv(\sigma)\bigr)$
follows from the observation that $\Phi\inv(\sigma)$ is equal to
$\Phi\inv(\bar\sigma)\setminus\bigcup_{\tau<\sigma}\Phi\inv(\tau)$.
\end{proof}

The moment map is determined up to an additive constant vector in
$\z^*$, where $\z$ is the Lie algebra of $Z$, the unit component of
the centre of $K$.  The choice of this constant does not affect the
imploded cross-section.  This is most easily seen from the direct-sum
decomposition $\k=\z\oplus[\k,\k]$, which, by identifying $\z^*$ with
the annihilator of $[\k,\k]$ and $[\k,\k]^*$ with the annihilator of
$\z$, gives rise to a decomposition $\k^*=\z^*\oplus[\k,\k]^*$.
Correspondingly, the Weyl chamber decomposes into a product of a
vector space and a proper cone,
$\t^*_+=\z^*\times(\t^*_+\cap[\k,\k]^*)$.  In fact, this argument
proves the first assertion of the following lemma.  The second
assertion is clear.

\begin{lemma}\label{lemma;semisimple}
The imploded cross-section of $M$ with respect to the $K$-action is
the same as the imploded cross-section with respect to the
$[K,K]$-action.  Likewise\upn, replacing $K$ by a finite cover does
not alter the imploded cross-section.
\end{lemma}

To obtain more detailed information we invoke the cross-section
theorem, which is essentially due to Guillemin and Sternberg; cf.\
\cite{guillemin-lerman-sternberg;symplectic-fibrations}.  The version
stated below incorporates some refinements made by Lerman et al.\
\cite{lerman;nonabelian-convexity}.  Consider a face $\sigma$ of
$\t^*_+$ and the $K_\sigma$-stable subset $\eu
S_\sigma=K_\sigma\cdot\star\sigma$ of $\k^*$, where $\star\sigma$
denotes the open star $\bigcup_{\tau\ge\sigma}\tau$ of $\sigma$.  It
is well-known that $\eu S_\sigma$ is a slice for the coadjoint action
(i.e.\ $K\eu S_\sigma$ is open and $K$-equivariantly diffeomorphic to
the associated bundle $K\times^{K_\sigma}\eu S_\sigma$), in fact the
largest possible slice containing all points of orbit type
$(K_\sigma)$.  The \emph{symplectic cross-section} of $M$ over
$\sigma$ is the subset
$$
M_\sigma=\Phi\inv(\eu S_\sigma).
$$
Note that $\Phi(M_\sigma)\subset\eu S_\sigma\subset\k_\sigma^*$ and
that the saturation $KM_\sigma$ of $M_\sigma$ is open in $M$.  The
cross-section $M_{\sigma\prin}$ is called \emph{principal}.

\begin{theorem}[symplectic cross-sections]\label{theorem;cross}
Let $\sigma$ be an open face of $\t^*_+$.
\begin{enumerate}
\item\label{item;cross}
The cross-section $M_\sigma$ is a connected $K_\sigma$-stable
symplectic submanifold of $M$.  The $K_\sigma$-action on $M_\sigma$ is
Hamiltonian with moment map $\Phi_\sigma=\Phi|_{M_\sigma}$.
\item\label{item;multiplication}
The multiplication map $K\times M_\sigma\to M$ induces a
symplectomorphism $K\times^{K_\sigma}M_\sigma\to KM_\sigma$.  If
$M_\sigma$ is nonempty\upn, then $KM_\sigma$ is dense in $M$.
\item\label{item;principal}
The commutator subgroup of $K_{\sigma\prin}$ acts trivially on
$M_{\sigma\prin}$.
\end{enumerate}
\end{theorem}

It follows from \eqref{item;principal} that if
$\Phi(m_1)\in\sigma\prin$, then $m_1\sim m_2$ implies $m_1=m_2$.  From
\eqref{item;cross} it then follows that
$\pi\bigl(\Phi\inv(\sigma\prin)\bigr)$ is connected and open and from
\eqref{item;multiplication} that it is dense.

\begin{corollary}\label{corollary;open-dense}
The restriction of $\pi$ to $\Phi\inv(\sigma\prin)$ is a homeomorphism
onto its image.  The image is connected\upn, and open and dense in
$M\impl$.  Hence $M\impl$ is connected.
\end{corollary}

By Theorem \ref{theorem;cross}\eqref{item;cross}, the composition of
$\Phi_\sigma$ with the projection
$\k_\sigma^*\to[\k_\sigma,\k_\sigma]^*$ is a moment map for the action
of $[K_\sigma,K_\sigma]$ on $M_\sigma$.  Its zero fibre is
$\Phi\inv(\z_\sigma^*)\cap M_\sigma=\Phi\inv(\z_\sigma^*\cap\eu
S_\sigma) =\Phi\inv(\sigma)$.  The decomposition
\eqref{equation;decomposition} can therefore be written more
insightfully as follows.

\begin{corollary}
The imploded cross-section partitions into symplectic quotients\upn,
each of which is locally closed\upn,
\begin{equation}\label{equation;quotients}
M\impl =\coprod_{\sigma\in\Sigma}M_\sigma\quot{[K_\sigma,K_\sigma]}.
\end{equation}
\end{corollary}

Here the notation $\quot[\lambda]$ indicates symplectic reduction at
level $\lambda$, the subscript being usually omitted when it is $0$.
For instance, the piece corresponding to the smallest face of the Weyl
chamber is $\Phi\inv(\z^*)/[K,K]=M\quot{[K,K]}$.

Not all the pieces of the partition \eqref{equation;quotients} are
symplectic manifolds, but a decomposition into symplectic manifolds
can be obtained by subdividing each of the pieces according to orbit
type.  For any $\sigma\in\Sigma$ and any closed subgroup $H$ of
$K'=[K_\sigma,K_\sigma]$, let
$$
M_{\sigma,(H)}=\{\,m\in M_\sigma\mid\text{$K'_m$ is conjugate within
$K'$ to $H$}\,\}
$$
be the stratum of orbit type $(H)$ in the $K'$-manifold $M_\sigma$.
Here $K'_m$ is the stabilizer of $m$ with respect to the $K'$-action.
By \cite[Theorem 2.1]{sjamaar-lerman;stratified}, the intersection
$\Phi\inv(\sigma)\cap M_{\sigma,(H)}$ is a smooth $K'$-stable
submanifold of $M_\sigma$ and the null-foliation of the symplectic
form restricted to this submanifold is exactly given by the
$K'$-orbits.  Hence the quotient
\begin{equation}\label{equation;stratum}
\bigl(\Phi\inv(\sigma)\cap M_{\sigma,(H)}\bigr)\big/K'
\end{equation}
is a symplectic manifold in a natural way.  It is more convenient to
work with the connected components of these manifolds instead.  Let
$\{\,X_i\mid i\in\ca I\,\}$ be the collection of components of all
manifolds of the form \eqref{equation;stratum}, where $\sigma$ ranges
over all faces of $\t^*_+$ and $(H)$ over all conjugacy classes of
subgroups of $[K_\sigma,K_\sigma]$.  We call the $X_i$ \emph{strata},
although it will not be proved until Section
\ref{section;stratification} that they form a stratification of
$M\impl$.  There is a partial order on the index set $\ca I$ defined
by $i\le j$ if $X_i\subset\bar X_j$.  By Corollary
\ref{corollary;open-dense}, $\ca I$ has a unique maximal element.
Moreover, the orbit type decomposition of any manifold with a smooth
action of a compact Lie group is locally finite.  Together with the
fact that the quotient map $\pi$ is proper (Lemma
\ref{lemma;hausdorff}), this implies that the collection of strata is
locally finite.  We have proved the following result.

\begin{theorem}\label{theorem;decomposition}
The imploded cross-section is a locally finite disjoint union of
locally closed connected subspaces\upn, each of which is a symplectic
manifold\upn,
\begin{equation}\label{equation;stratification}
M\impl=\coprod_{i\in\ca I}X_i.
\end{equation}
There is a unique open stratum\upn, which is dense in $M\impl$ and
symplectomorphic to the principal cross-section of $M$.
\end{theorem}

\section{Abelianization}\label{section;abelian}

The imploded cross-section of a Hamiltonian $K$-manifold can be viewed
as its abelianization in a sense which we shall now make precise.

Let $\ca X$ be a topological space with a decomposition $\ca
X=\coprod_{i\in\ca I}\ca X_i$ into connected subspaces $\ca X_i$, each
of which is equipped with the structure of a smooth manifold and a
symplectic form $\omega_i$.  A continuous action of a Lie group $\ca
G$ on $\ca X$ is \emph{Hamiltonian} if it preserves the decomposition
and is smooth on each $\ca X_i$, and if there exists a continuous
$\Ad^*$-equivariant map $\Phi_{\ca X}\colon\ca X\to(\Lie\ca G)^*$, the
\emph{moment map}, such that $\Phi_{\ca X}|_{\ca X_i}$ is a moment map
in the usual sense for the $\ca G$-action on $\ca X_i$ for all
$i\in\ca I$.  The triple $\bigl(\ca X,\{\,(\ca X_i,\omega_i)\mid i\in
\ca I\,\},\Phi_{\ca X}\bigr)$ is a \emph{Hamiltonian $\ca G$-space}.
An \emph{isomorphism} from $\ca X$ to a second Hamiltonian $\ca
G$-space $\bigl(\ca Y,\{\,(\ca Y_j,\omega_j)\mid j\in\ca
J\,\},\Phi_{\ca Y}\bigr)$ is a pair $(F,f)$, where $F$ is a
homeomorphism $\ca X\to\ca Y$ and $f$ is a bijection $\ca I\to\ca J$
subject to the following conditions: $F$ is equivariant, $\Phi_{\ca
X}=\Phi_{\ca Y}\circ F$, $F$ maps $\ca X_i$ symplectomorphically onto
$\ca Y_{f(i)}$ for all $i\in\ca I$.

\begin{example}\label{example;affine}
Let $V$ be a finite-dimensional unitary $K$-module with inner product
$\langle\cdot,\cdot\rangle$.  This is a Hamiltonian $K$-manifold with
symplectic form $\omega_V$ and moment map $\Phi_V$ given by
\begin{equation}\label{equation;quadratic}
\omega_V=-\Im\langle\cdot,\cdot\rangle\qquad
\text{and}\qquad\Phi_V^\xi(v)=\frac1{2}\,\omega_V(\xi v,v),
\end{equation}
respectively.  Let $\ca X$ be a $K$-stable irreducible complex
algebraic subvariety of $V$, regarded as a subspace for the classical
topology on $V$.  There is a natural minimal decomposition of $\ca X$
into nonsingular $K$-stable algebraic subvarieties, each of which
inherits a symplectic form and a moment map from the ambient space
$V$.  Let us call the Hamiltonian $K$-space $\ca X$ thus obtained an
\emph{affine Hamiltonian $K$-space}.  The topology on $\ca X$ and its
decomposition into manifolds depend only on its coordinate ring $\eu
A$, and the $K$-action is determined by the $K$-module structure of
$\eu A$.  The symplectic forms and the moment map, however, depend on
the embedding into $V$ and the inner product on $V$.  When the
embedding and the inner product are fixed, we will sometimes abuse
notation and write $\ca X=\Spec\eu A$ to indicate that $\ca X$ is the
affine Hamiltonian $K$-space whose underlying variety is the
subvariety of $V$ with coordinate ring $\eu A$.
\end{example}

For $\gamma\in(\Lie\ca G)^*$, the \emph{symplectic quotient} (or
\emph{reduced space}) at level $\gamma$ is the topological space $\ca
X\quot[\gamma]\ca G =\Phi_{\ca X}\inv(\gamma)/\ca G_\gamma$.  If $\ca
G$ is compact (or if $\ca G_\gamma$ acts properly on $\Phi_{\ca
X}\inv(\gamma)$), the symplectic quotient can be decomposed into
connected smooth symplectic manifolds by partitioning each of the
pieces $\ca X_i$ according to $\ca G$-orbit type, $\ca
X_i=\coprod_{(\ca H)}\ca X_{i,(\ca H)}$, forming the symplectic
manifolds $\bigl(\Phi_{\ca X}\inv(\gamma)\cap\ca X_{i,(\ca
H)}\bigr)\big/\ca G_\gamma$ as in \eqref{equation;stratum}, and
subdividing these into their connected components.

\begin{example}\label{example;projective}
Let $V$ be as in Example \ref{example;affine} and let $\P(V)$ be the
associated projective space.  As a symplectic manifold $\P(V)$ is
isomorphic to the quotient $V\quot[-1]S^1$, where $S^1$ acts by scalar
multiplication.  Let $\ca X$ be a $K$-stable irreducible complex
algebraic subvariety of $\P(V)$.  By partitioning into a minimal set
of nonsingular algebraic subvarieties as in Example
\ref{example;affine}, we see that $\ca X$ is a Hamiltonian $K$-space,
called a \emph{projective Hamiltonian $K$-space}.  
The affine cone $\tilde{\ca X}=\Spec\eu A\subset V$ is an affine
Hamiltonian $K$-space, and we have $\ca X\cong\tilde{\ca
X}\quot[-1]S^1$ as Hamiltonian $K$-spaces.
\end{example}

Now let $(M,\omega,K)$ be an arbitrary Hamiltonian $K$-manifold.  We
claim that the imploded cross-section $M\impl$, equipped with the
decomposition \eqref{equation;stratification}, is a Hamiltonian
$T$-space.  Indeed, the preimage $\Phi\inv(\t^*_+)$ is stable under
the action of the maximal torus.  In addition, $m_1\sim m_2$ implies
$tm_1\sim tm_2$ for all $t\in T$ and $m_1$, $m_2\in\Phi\inv(\t^*_+)$,
because $T$ normalizes each of the groups $[K_\sigma,K_\sigma]$.  Thus
the action of $T$ descends to a continuous action on $M\impl$.
This action is Hamiltonian with moment map $\Phi\impl$.  The easiest
way to see this is to use the following alternative definition of the
$T$-action: for each $\sigma$ the $K_\sigma$-action on $M_\sigma$
descends to a $K_\sigma/[K_\sigma,K_\sigma]$-action on the quotient
$M_\sigma\quot{[K_\sigma,K_\sigma]}$.  Via the canonical surjective
map $T\to K_\sigma/[K_\sigma,K_\sigma]$ this induces a (non-effective)
$T$-action on $M_\sigma\quot{[K_\sigma,K_\sigma]}$.  Again because $T$
normalizes $[K_\sigma,K_\sigma]$, this action preserves the
$[K_\sigma,K_\sigma]$-orbit type strata and is Hamiltonian on each
such stratum.  The moment maps are induced by the restrictions of
$\Phi$ to the manifolds $\Phi\inv(\sigma)\cap M_{\sigma,(H)}$; in
other words they are the restrictions of $\Phi\impl$ to the various
strata.

Symplectic reduction of $M\impl$ with respect to $T$ turns out to be
the same as symplectic reduction of $M$ with respect to $K$.  Namely,
let $\lambda\in\t^*_+$ and let $\sigma$ be the face of $\t^*_+$
containing $\lambda$.  Then
$\Phi\impl\inv(\lambda)=\Phi\inv(\lambda)/[K_\sigma,K_\sigma]$, so
there is a quotient map $\Phi\inv(\lambda)\to\Phi\impl\inv(\lambda)$.

\begin{theorem}\label{theorem;abelian}
For every $\lambda\in\t^*_+$\upn, the canonical map
$\Phi\inv(\lambda)\to\Phi\impl\inv(\lambda)$ induces an isomorphism of
symplectic quotients $M\quot[\lambda]K\cong M\impl\quot[\lambda]T$.
\end{theorem}

\begin{proof}
Let $\sigma$ be the face containing $\lambda$.  Assume first that all
points in $\Phi\inv(\lambda)$ are of the same orbit type for the
action of $K_\lambda=K_\sigma$, so that $M\quot[\lambda]K$ is a smooth
symplectic manifold.  By reduction in stages, it is naturally
symplectomorphic to the iterated quotient
\begin{equation}\label{equation;iterate}
\bigl(M_\sigma\quot[0][K_\sigma,K_\sigma]\bigr)\bigquot[\lambda]T.
\end{equation}
Since $\Phi\impl\inv(\lambda)\subset
M_\sigma\quot{[K_\sigma,K_\sigma]}$ and the restriction of $\Phi\impl$
to $M_\sigma\quot{[K_\sigma,K_\sigma]}$ is the moment map for the
$T$-action, it is clear that \eqref{equation;iterate} is naturally
symplectomorphic to $M\impl\quot[\lambda]T$.  If $\Phi\inv(\lambda)$
consists of more than one stratum, the same argument, using stratified
reduction in stages (see \cite[\S4]{sjamaar-lerman;stratified}), shows
that the quotient map $\Phi\inv(\lambda)\to\Phi\impl\inv(\lambda)$
induces a homeomorphism $M\quot[\lambda]K\to M\impl\quot[\lambda]T$,
which maps strata symplectically onto strata.
\end{proof}

\begin{example}\label{example;universal-torus}
For all $\lambda\in\t^*_+$, $T^*K\quot[\lambda]K\cong K\lambda$, the
coadjoint orbit through $\lambda$, so $(T^*K)\impl\quot[\lambda]T\cong
K\lambda$.
\end{example}

\section{The universal imploded cross-section}\label{section;universal}

As an example we determine the cross-sections and the imploded
cross-section of the cotangent bundle $T^*K$.  It turns out that the
space $(T^*K)\impl$ has a universal property, which greatly
facilitates calculations involving symplectic implosion.  Another
interesting feature is its homogeneous structure.  We shall see that
if $K$ is semisimple, $(T^*K)\impl$ is a cone over a compact space,
which stratifies into contact manifolds.

Consider the actions of $K$ on itself given by $\ca L_gk=gk$ and $\ca
R_gk=kg\inv$, which both lift to a Hamiltonian actions on $T^*K$.
Identify $T^*K$ with $K\times\k^*$ by means of left translations.
Then the actions are given by $\ca L_g(k,\lambda)=(gk,\lambda)$ and
$\ca R_g(k,\lambda)=(kg\inv,g\lambda)$, where $g\lambda$ is an
abbreviation for $\Ad^*(g)(\lambda)$.  The moment maps (with respect
to the symplectic form $\omega=d\beta$, where $\beta$ is the canonical
one-form) are respectively
$$
\Phi_{\ca L}(k,\lambda)=-k\lambda,\qquad\Phi_{\ca R}(k,\lambda)
=\lambda.
$$
The inversion map $k\mapsto k\inv$ intertwines the left and right
actions.  Its cotangent lift, given by
$(k,\lambda)\mapsto(k\inv,-k\lambda)$, is a symplectic involution of
$T^*K$, which intertwines $\Phi_{\ca L}$ and $\Phi_{\ca R}$.
Therefore the cross-sections for the left action are canonically
isomorphic to those for the right action.  For simplicity let us use
$\Phi_{\ca R}$.  Clearly
\begin{gather}
\label{equation;cross}
\Phi_{\ca R}\inv(\eu S_\sigma)=K\times\eu S_\sigma,\\
\label{equation;implode}
(T^*K)\impl =\coprod_{\sigma\in\Sigma}(K\times\eu
S_\sigma)\quot{[K_\sigma,K_\sigma]}
=\coprod_{\sigma\in\Sigma}\frac{K}{[K_\sigma,K_\sigma]}\times\sigma,
\end{gather}
so in this example the decompositions \eqref{equation;quotients} and
\eqref{equation;stratification} are identical.

As stated in Theorem \ref{theorem;cross}\eqref{item;cross}, the
cross-section \eqref{equation;cross} inherits a symplectic structure
from $T^*K$.  Here is an alternative construction of the symplectic
form.  For each face $\sigma$ the subalgebra $\k_\sigma$ is equal to
the subspace of $Z_\sigma$-invariants in $\k$, where $Z_\sigma$ is the
identity component of the centre of $K_\sigma$.  We therefore have a
canonical $K_\sigma$-invariant projection $\k\to\k_\sigma$ and hence a
connection $\theta$ on the principal $K_\sigma$-bundle
\begin{equation}\label{equation;principal}
K_\sigma\to K\to K/K_\sigma.
\end{equation}
A connection $\vartheta\in\Omega^1(\ca P,\Lie\ca G)$ on a principal
bundle $\ca P$ for a Lie group $\ca G$ is \emph{fat} at
$\gamma\in(\Lie\ca G)^*$ if the two-form
$\langle\gamma,d\vartheta\rangle$ is nondegenerate on the horizontal
subspaces of $\ca P$.  This is a necessary and sufficient condition
for the closed two-form $d\langle\pr_2,\vartheta\rangle$ on $\ca
P\times(\Lie\ca G)^*$ to be nondegenerate at $P\times\{\gamma\}$,
where $\pr_2\colon\ca P\times(\Lie\ca G)^*\to(\Lie\ca G)^*$ is the
projection.

\begin{lemma}%
[{\cite[Corollary
2.3.8]{guillemin-lerman-sternberg;symplectic-fibrations}}]
\label{lemma;fat}
The canonical connection $\theta$ on the bundle
\eqref{equation;principal} is fat at $\lambda\in\k_\sigma^*$ if and
only if $\lambda\in\eu S_\sigma$.
\end{lemma}

Therefore the form $d\langle\pr_2,\theta\rangle$ on
$K\times\k_\sigma^*$ is symplectic on $K\times\eu S_\sigma$.  It is
straightforward to check that $\langle\pr_2,\theta\rangle$ is equal to
the restriction to $K\times\eu S_\sigma$ of the canonical one-form
$\beta$.  Hence $d\langle\pr_2,\theta\rangle$ is equal to the
restriction of $d\beta=\omega$.

The symplectic form on the stratum $(K\times\eu
S_\sigma)\quot{[K_\sigma,K_\sigma]}$ of $(T^*K)\impl$ can therefore be
interpreted as the form obtained by reducing $(K\times\eu
S_\sigma,d\langle\pr_2,\theta\rangle)$ with respect to
$[K_\sigma,K_\sigma]$.  A third alternative is to note that the
projection $\k\to\k_\sigma$ descends to a $Z_\sigma$-equivariant
projection $\k/[\k_\sigma,\k_\sigma]\to\k_\sigma/[\k_\sigma,\k_\sigma]
\cong\z_\sigma$, whence we obtain a connection
$\bar\theta\in\Omega^1\bigl(K/[K_\sigma,K_\sigma],\z_\sigma\bigr)$ on
the torus bundle
\begin{equation}\label{equation;subprincipal}
K_\sigma/[K_\sigma,K_\sigma]\to K/[K_\sigma,K_\sigma]\to K/K_\sigma.
\end{equation}
Put $\beta_\sigma =\langle\pr_2,\bar\theta\rangle$ and
$\omega_\sigma=d\langle\pr_2,\bar\theta\rangle$, where $\pr_2$ now
stands for the projection
$K/[K_\sigma,K_\sigma]\times\sigma\to\sigma$.  Lemma \ref{lemma;fat}
implies that $\bar\theta$ is fat at all points of $\eu
S_\sigma\cap\z_\sigma^*=\sigma$, and therefore $\omega_\sigma$ is a
symplectic form on $K/[K_\sigma,K_\sigma]\times\sigma=(K\times\eu
S_\sigma)\quot{[K_\sigma,K_\sigma]}$.  The following result asserts
that this form is the same as the one defined by symplectic implosion.

\begin{lemma}\label{lemma;exact}
Let $p$ be the orbit map $K\times\sigma\to(K\times\eu
S_\sigma)\quot{[K_\sigma,K_\sigma]}$.  Then
$\beta|_{K\times\sigma}=p^*\beta_\sigma$ and
$\omega|_{K\times\sigma}=p^*\omega_\sigma$.
\end{lemma}

\begin{proof}
The second equality is immediate from the first.  Because of
$K$-invariance, the first equality need only be checked at points of
the form $(1,\lambda)\in K\times\sigma$.  Let $(\xi,\mu)\in
T_{(1,\lambda)}(K\times\sigma)=\k\times\z^*_\sigma$.  Then
$\beta_{(1,\lambda)}(\xi,\mu)=\lambda(\xi)$ by the definition of
$\beta$.  (We identify $\z^*_\sigma$ with the annihilator of
$[\k_\sigma,\k_\sigma]\cap\t$ in $\t^*$, so that $\lambda(\xi)$ is
well-defined.)  Moreover,
$$
(p^*\beta_\sigma)_{(1,\lambda)}(\xi,\mu)
=(\beta_\sigma)_{p(1,\lambda)}(\xi\bmod[\k_\sigma,\k_\sigma],\mu)
=\lambda\bigl(\bar\theta(\xi\bmod[\k_\sigma,\k_\sigma])\bigr),
$$
which is equal to $\lambda(\xi)$ because $\lambda\in\z_\sigma^*$ and
$\bar\theta(\xi\bmod[\k_\sigma,\k_\sigma])$ is equal to the projection
of $\xi$ onto $\z_\sigma$.
\end{proof}

Not only is $(T^*K)\impl$ a Hamiltonian $T$-space for the $T$-action
induced by the right $K$-action, but the left $K$-action on $T^*K$
descends to a $K$-action on $(T^*K)\impl$ given by $g\pi_{\ca
R}(k,\lambda)=\pi_{\ca R}(gk,\lambda)$, which is Hamiltonian as well.
(Here $\pi_{\ca R}$ denotes the quotient map $\Phi_{\ca R}\inv(\t^*_+)
\to(T^*K)\impl$.)  Its moment map is induced by $\Phi_{\ca L}$ and is
for simplicity also denoted by $\Phi_{\ca L}$.  Clearly the two
actions commute, so that $(T^*K)\impl$ is a Hamiltonian $K\times
T$-space with moment map $\Phi_{\ca L}\times(\Phi_{\ca R})\impl$.

\begin{example}\label{example;su2}
Let $K=\SU(2)$, which we shall identify with $S^3\subset\H$, the unit
quaternions.  The Lie algebra is then the set of imaginary
quaternions, the unit circle $S^1\subset\C$ is a maximal torus, and
the fibration \eqref{equation;principal} is none other than the Hopf
fibration $S^1\to S^3\to S^2$.  The symplectic form on
$K\times(\t^*_+)^\circ=S^3\times(0,\infty)$ is $d(t\theta)$, with $t$
being the standard coordinate on $(0,\infty)$.  The map
$(z,t)\mapsto\sqrt{2t}\,z$ is a symplectomorphism from
$S^3\times(0,\infty)$ onto $\H\setminus\{0\}$ with its standard
symplectic structure.  The remaining stratum, corresponding to
$\sigma=\{0\}$, consists of a single point.  Consequently, the
continuous map $F\colon T^*K\to\H=\C^2$ defined by
$F(k,\lambda)=\sqrt{2\,\norm{\lambda}}\,k$ induces a continuous
bijection $(T^*K)\impl\to\C^2$, which is a homeomorphism because
$(T^*K)\impl$ is locally compact Hausdorff.  In this sense, the
isolated singularity is removable and the imploded cross-section is a
symplectic $\C^2$.  Modulo this identification, the left $K$-action on
$(T^*K)\impl$ is the standard representation of $\SU(2)$ on $\C^2$ and
the right $T$-action is given by $t\cdot z=t\inv z$.  This example
will be generalized to arbitrary $K$ in Section \ref{section;kaehler}.
\end{example}

Now let $(M,\omega,\Phi)$ be any Hamiltonian $K$-manifold and define
$j\colon M\to M\times T^*K$ by $j(m)=(m,1,\Phi(m))$.

\begin{lemma}\label{lemma;universal}
\begin{enumerate}
\item\label{item;universal-centre}
The map $j$ is a symplectic embedding and induces an isomorphism of
Hamiltonian $K$-manifolds\upn,
$$
\bar\jmath\colon M\to(M\times T^*K)\quot K.
$$
Here the right-hand side is the quotient with respect to the diagonal
$K$-action\upn, where $K$ acts on the left on $T^*K$.  The $K$-action
on $(M\times T^*K)\quot K$ is the one induced by the right action on
$T^*K$.
\item\label{item;universal}
For every face $\sigma$\upn, $j$ maps $M_\sigma$ into $M\times
K\times\eu S_\sigma$ and induces an isomorphism of Hamiltonian
$K_\sigma$-manifolds\upn,
$$
\bar\jmath_\sigma\colon M_\sigma\to(M\times K\times\eu S_\sigma)\quot
K,
$$
where the quotient is taken as in \eqref{item;universal-centre}.
\end{enumerate}
\end{lemma}

\begin{proof}
The map $m\mapsto(1,\Phi(m))$ sends $M$ to the Lagrangian $\k^*\subset
T^*K$, so $j$ is a symplectic embedding.  The moment map for the
diagonal $K$-action on $M\times T^*K$ is given by
$\Psi(m,k,\lambda)=\Phi(m)-k\lambda$, so $j$ maps $M$ into
$\Psi\inv(0)$.  One checks readily that the induced map $\bar\jmath$
is a diffeomorphism.  It is a symplectomorphism because $j$ is a
symplectic embedding.  Moreover, $j(km)=(km,1,k\Phi(m)) =\ca L_k\ca
R_kj(m)$ and $\Phi_{\ca R}(j(m))=\Phi(m)$, so $\bar\jmath$ is
$K$-equivariant and intertwines the $K$-moment maps on $M$ and
$(M\times T^*K)\quot K$.  This proves \eqref{item;universal-centre}.
\eqref{item;universal} is proved in exactly the same way.
\end{proof}

Observe that $j$ maps $\cross$ into $\coprod_\sigma M\times
K\times\sigma$ and therefore induces a continuous map $\cross\to
M\times(T^*K)\impl$. Because of the following result we call
$(T^*K)\impl$ the \emph{universal imploded cross-section}.  See
Section \ref{section;abelian} for the definition of an isomorphism of
Hamiltonian $T$-spaces.

\begin{theorem}\label{theorem;universal}
The map $j$ induces an isomorphism of Hamiltonian $T$-spaces\upn,
$$
\bar\jmath\impl\colon M\impl\to\bigl(M\times(T^*K)\impl\bigr)\bigquot
K,
$$
where the quotient is taken with respect to the diagonal $K$-action.
\end{theorem}

\begin{proof}
The $K$-moment map on $M\times(T^*K)\impl$ is given by
$\Psi(m,\pi_{\ca R}(k,\lambda))=\Phi(m)-k\lambda$, so $j$ maps
$\cross$ into $\Psi\inv(0)$.  Moreover, $m_1\sim m_2$ implies
$j(m_1)=j(m_2)$, so $j$ induces a continuous map
$M\impl\to\Psi\inv(0)$.  Upon quotienting by $K$ we obtain the map
$\bar\jmath\impl$.  As in the proof of Lemma \ref{lemma;universal} one
checks that $\bar\jmath\impl$ is a homeomorphism which is
$T$-equivariant and intertwines the $T$-moment maps on $M\impl$ and
$\bigl(M\times(T^*K)\impl\bigr)\bigquot K$.  Note that
$\bar\jmath\impl$ restricts to a map
\begin{equation}\label{equation;induced}
M_\sigma\quot{[K_\sigma,K_\sigma]}\to\bigl(M\times(K\times\eu
S_\sigma)\quot{[K_\sigma,K_\sigma]}\bigr)\bigquot K.
\end{equation}
This is none other than the map induced, upon reduction with respect
to $[K_\sigma,K_\sigma]$, by the map $\bar\jmath_\sigma$ defined in
Lemma \ref{lemma;universal}\eqref{item;universal}.  Because
$\bar\jmath_\sigma$ is an isomorphism of Hamiltonian
$K_\sigma$-manifolds, it preserves the $[K_\sigma,K_\sigma]$-orbit
types and therefore \eqref{equation;induced} maps strata onto strata
and is a symplectomorphism on each stratum.
\end{proof}

\begin{remark}\label{remark;universal}
Let $S_\sigma=\Phi_{\ca L}\inv(\sigma)$ denote the stratum of
$(T^*K)\impl$ corresponding to a face $\sigma\in\Sigma$.  Furthermore,
let $\tau=\sigma\prin$ be the principal face of $M$.  Then the closure
of $S_\tau$ is equal to $\coprod_{\sigma\le\tau}S_\sigma$.  Since the
moment polytope of $M$ is contained in $\bar\tau
=\coprod_{\sigma\le\tau}\sigma$, the proof of Theorem
\ref{theorem;universal} shows that $j$ induces an isomorphism
$M\impl\cong(M\times\bar S_\tau)\quot K$.
\end{remark}

\begin{example}\label{example;multiplicity-free}
Let $M$ be a point.  Then the theorem asserts that
$(T^*K)\impl\quot[\lambda]K$ consists of a single point for all
$\lambda\in\k^*$.  In particular, $(T^*K)\impl$ is a multiplicity-free
$K$-space.
\end{example}

As an application, consider a Hamiltonian action of a second compact
Lie group $H$ on $M$, which commutes with $K$.  The $H$-moment map on
$M$ is $K$-invariant and therefore induces a continuous map on
$M\impl$ (implosion with respect to the $K$-action).  This map is a
moment map for the $H$-action induced on $M\impl$.  Using Theorem
\ref{theorem;universal} and reduction in stages we conclude that
reduction commutes with implosion.

\begin{corollary}\label{corollary;reduction-implosion}
The Hamiltonian $T$-spaces $(M\quot[\eta]H)\impl$ and
$M\impl\quot[\eta]H$ are isomorphic for every $\eta\in\lie l^*$.
\end{corollary}

\begin{example}\label{example;universal}
Consider the $K\times K$-space $M\times T^*K$, where the first copy of
$K$ acts diagonally (and by left multiplication on $T^*K$), and the
second copy acts by right multiplication on $T^*K$.  Reducing with
respect to the first copy, imploding with respect to the second, and
applying Lemma \ref{lemma;universal}\eqref{item;universal-centre} we
find
$$
M\impl\cong\bigl((M\times T^*K)\quot K\bigr)\impl
\cong\bigl(M\times(T^*K)\impl\bigr)\bigquot K.
$$
In other words, Corollary \ref{corollary;reduction-implosion} is
equivalent to Theorem \ref{theorem;universal}.
\end{example}

\begin{example}\label{example;semisimple}
If $H$ is finite, then $(M/H)\impl\cong M\impl/H$.  In particular, let
$\Gamma$ be the finite central subgroup $Z\cap[K,K]$ of $K$.  Then
$K=(Z\times[K,K])/\Gamma$ and $T^*K=T^*Z\times^\Gamma T^*[K,K]$.
Since implosion relative to $K$ is the same as implosion relative to
$[K,K]$, we see that
\begin{equation}\label{equation;equisingular}
(T^*K)\impl\cong(T^*Z\times T^*[K,K])\impl/\Gamma\cong
T^*Z\times^\Gamma(T^*[K,K])\impl,
\end{equation}
a bundle with fibre $(T^*[K,K])\impl$ over the cotangent bundle of the
torus $Z/\Gamma$.  In turn, $(T^*[K,K])\impl$ can be written as
$(T^*[K,K]\sptilde)\impl/\Delta$, where $\Delta$ is the fundamental
group of $[K,K]$ and $[K,K]\sptilde$ its universal cover.  For
instance, if $K=\SO(3)=\SU(2)/\{\pm1\}$, then Example
\ref{example;su2} shows that $(T^*K)\impl$ is the symplectic orbifold
$\C^2/\{\pm\id\}$, and if $K=\U(2)$ we find
$(T^*K)\impl=\C^\times\times^{\{\pm\id\}}\C^2$.
\end{example}

Because of \eqref{equation;equisingular}, to describe the
singularities of $(T^*K)\impl$ we can focus our attention on the fibre
$(T^*[K,K])\impl$.  So let us assume for the remainder of this section
that $K$ is semisimple.  For $t\ge0$ let $\ca A_t$ be fibrewise
multiplication by $t$ in $T^*K$.  Let
$$
T^\times K=T^*K\setminus\{\text{zero section}\}
=K\times(\k^*\setminus\{0\})
$$
be the punctured cotangent bundle and $\R^\times=(0,\infty)$ the
multiplicative group of positive reals.  Then for $t>0$ $\ca A_t$
defines a proper free action of $\R^\times$ on $T^\times K$.  Let
$\zeta$ be its infinitesimal generator.  Then $\zeta$ is a Liouville
vector field, i.e.\ $\ca L(\zeta)\omega=\omega$, and therefore we call
$\ca A$ a \emph{Liouville action}.  We can write the canonical
one-form as $\beta=\iota(\zeta)\omega$.  Recall that if $\nu$ is a
contact form on a manifold, then the associated \emph{Reeb vector
field} is the vector field $\chi$ defined by $\iota(\chi)d\nu=0$ and
$\iota(\chi)\nu=1$.

\begin{lemma}\label{lemma;contact}
Let $\zeta$ be a global Liouville vector field on a symplectic
manifold $(M,\omega)$ and let $\beta$ be the potential one-form
$\iota(\zeta)\omega$.  Let $\Xi$ be a symplectic vector field on $M$
that commutes with $\zeta$.
\begin{enumerate}
\item\label{item;hamiltonian}
The function $\varphi=-\iota(\Xi)\beta$ is a Hamiltonian for $\Xi$ and
satisfies $\ca L(\zeta)\varphi=\varphi$.
\item\label{item;contact}
Any $c\ne0$ is a regular value of $\varphi$\upn, the restriction of
$\beta$ to the hypersurface $\varphi\inv(c)$ is a contact form\upn,
and the restriction of $-c\inv\Xi$ is its Reeb vector field.
\end{enumerate}
\end{lemma}

\begin{proof}
By assumption
$$
\ca L(\Xi)\beta=\ca L(\Xi)\iota(\zeta)\omega=\bigl(\iota(\zeta)\ca
L(\Xi)+\iota([\Xi,\zeta])\bigr)\omega=0,
$$
and hence
\begin{gather*}
d\varphi=-d\iota(\Xi)\beta=-\ca L(\Xi)\beta+\iota(\Xi)d\beta
=\iota(\Xi)d\beta=\iota(\Xi)\omega,\\
\ca L(\zeta)\varphi=\iota(\zeta)d\varphi=\iota(\zeta)\iota(\Xi)\omega
=-\iota(\Xi)\beta=\varphi,
\end{gather*}
so \eqref{item;hamiltonian} holds.  It follows that
$d\varphi_m(\zeta)=\ca L(\zeta)\varphi_m=\varphi(m)$, so any $c\ne0$
is a regular value and $\zeta$ is transverse to $\varphi\inv(c)$.  It
is now easy to see that $\beta$ is a contact form on $\varphi\inv(c)$;
see e.g.\ McDuff and Salamon \cite[Proposition
3.57]{mcduff-salamon;introduction}.  Furthermore, on $\varphi\inv(c)$
we have $\iota(\Xi)d\beta =\iota(\Xi)\omega=d\varphi=0$ and
$\iota(\Xi)\beta =-\varphi(m)=-c$, so $-c\inv\Xi|_{\varphi\inv(c)}$ is
the Reeb vector field, which proves \eqref{item;contact}.
\end{proof}

\begin{example}\label{example;cn}
Let $M=\C^n$ with its standard symplectic form and let
$\zeta=\frac1{2}\sum_i(x_i\partial/\partial x_i+y_i\partial/\partial
y_i)$ and $\Xi =\sum_i(-y_i\partial/\partial x_i+x_i\partial/\partial
y_i)$.  Then $\zeta$ is a Liouville vector field and generates the
action $(t,z)\mapsto\sqrt{t}\,z$, and $\Xi$ generates the standard
circle action.  Moreover, $\varphi(z)=-\frac1{2}\norm{z}^2$, so for
$c<0$ the hypersurfaces $\varphi\inv(c)$ are spheres and the orbits of
the Reeb vector field on $\varphi\inv(c)$ are the fibres of the Hopf
fibration.
\end{example}

\begin{example}\label{example;cotangent}
Let $M=T^\times K$ and choose $\Xi\in\t$.  Let $T$ act on $M$ from the
right and consider the vector field on $M$ induced by $\Xi$, which for
brevity we will also denote by $\Xi$.  Then
$\varphi(k,\lambda)=\lambda(\Xi)$ and $\Xi$ commutes with $\zeta=d\ca
A_t/dt|_{t=0}$.  Therefore
$\varphi\inv(-1)=K\times\{\,\lambda\mid\lambda(\Xi)=-1\,\}$ is a
hypersurface of contact type with Reeb vector field $\Xi$.
\end{example}

In Example \ref{example;cn} the hypersurface $\varphi\inv(-1)$ is
compact and $M$ is topologically a cone over $\varphi\inv(-1)$.  This
is obviously not the case in Example \ref{example;cotangent}, but we
shall now show that this can be remedied by imploding $M$.  Since $K$
is semisimple, $\Phi_{\ca R}\inv(\z^*)=\Phi_{\ca R}\inv(0)$ is the
zero section of $T^*K$, and therefore
$$
(T^\times K)\impl=(T^*K)\impl\setminus\{*\}
=\coprod_{\sigma\in\Sigma\setminus\{0\}}\!\!\!\!F_\sigma,
$$
where $F_\sigma=(K\times\eu S_\sigma)\quot{[K_\sigma,K_\sigma]}$ and
$\{*\}=F_{\{0\}}$ is the \emph{vertex} of $(T^*K)\impl$.  Because
$\Phi_{\ca R}$ is homogeneous and equivariant, the action $\ca A$
descends to an action on $(T^*K)\impl$, denoted also by $\ca A$, which
off the vertex is proper and free, and on each stratum $F_\sigma$ is a
Liouville action.  For each $F_\sigma$, let $\omega_\sigma$ be the
symplectic form, $\zeta_\sigma$ the Liouville vector field,
$\beta_\sigma$ the one-form $\iota(\zeta_\sigma)\omega_\sigma$ and
$\Xi_\sigma$ the vector field induced by $\Xi\in\t$.  The Hamiltonian
$\varphi=-\iota(\Xi)\beta$ descends to a continuous function on
$(T^*K)\impl$, also denoted by $\varphi$.  The functions
$\varphi_\sigma=\varphi|_{F_\sigma}$ satisfy
$\varphi_\sigma=-\iota(\Xi_\sigma)\beta_\sigma$.  The subset
$\varphi\inv(-1)$ of $(T^*K)\impl$ is called the \emph{link} of the
vertex and denoted by $\lk(*)$.  The infinite cone
$\bigl(X\times[\,0,\infty)\bigr)\big/\bigl(X\times\{0\}\bigr)$ over a
space $X$ is denoted by $C^\circ(X)$.

\begin{proposition}\label{proposition;homogeneous}
Assume that $K$ is semisimple and that $\Xi$ is in the interior of the
cone spanned by the negative coroots.  Then
\begin{enumerate}
\item\label{item;proper}
$\varphi$ is proper on $(T^*K)\impl$\upn, $\lk(*)$ is compact\upn,
$\varphi\le0$\upn, and $\varphi\inv(0)=\{*\}$\upn;
\item\label{item;subcontact}
for all $\sigma\ne0$\upn, the intersection $\lk(*)_\sigma =\lk(*)\cap
F_\sigma$ is a smooth manifold\upn, $\beta_\sigma$ restricts to a
contact form\upn, and $\Xi_\sigma$ to the Reeb vector field on
$\lk(*)_\sigma$\upn;
\item\label{item;section}
The link of $*$ is a global section of the principal
$\R^\times$-bundle $(T^\times K)\impl$.  The map
$f\colon\lk(*)\times\R^\times\to(T^\times K)\impl$ given by
$f(\pi_{\ca R}(k,\lambda),t)=\ca A_t(\pi_{\ca R}(k,\lambda))=\pi_{\ca
R}(k,t\lambda)$ is a stratum-preserving homeomorphism\upn, which on
every stratum is a diffeomorphism\upn, and satisfies
$f^*\omega_\sigma=d(t(\beta_\sigma|_{\lk(*)_\sigma}))$\upn;
\item\label{item;cone}
$f$ extends uniquely to a homeomorphism
$C^\circ(\lk(*))\to(T^*K)\impl$.
\end{enumerate}
\end{proposition}

\begin{proof}
The cone spanned by the \emph{positive} coroots is the dual of the
cone $\t^*_+$.  Since $-\Xi$ is in its interior and $K$ is semisimple,
the linear function $\lambda\mapsto\lambda(\Xi)$ is proper on
$\t^*_+$.  It follows that $\varphi(\pi_{\ca
R}(k,\lambda))=\lambda(\Xi)$ is proper on $(T^*K)\impl$.  The other
assertions in \eqref{item;proper} are now obvious.  To prove
\eqref{item;subcontact}, apply Lemma \ref{lemma;contact} to each
stratum $F_\sigma$.  \eqref{item;section} follows from the observation
that the simplex $\{\,\lambda\in\t^*_+\mid\lambda(\Xi)=-1\,\}$ is a
global section of the $\R^\times$-action $\lambda\mapsto t\lambda$ on
the punctured Weyl chamber $\t^*_+\setminus\{0\}$.  The equality
$f^*\omega_\sigma=d(t(\beta_\sigma|_{\lk(*)_\sigma}))$ is readily
checked on the tangent space $T_mF_\sigma$ at any $m\in\lk(*)_\sigma$
and therefore holds globally by homogeneity.  To prove
\eqref{item;cone}, observe that $(\pi_{\ca
R}(k,\lambda),t)\mapsto\pi_{\ca R}(k,t\lambda)$ defines a map
$\lk(*)\times[\,0,\infty)\to(T^*K)\impl$ which sends
$\lk(*)\times\{0\}$ to $*$ and therefore descends to a homeomorphism
$C^\circ(\lk(*))\to(T^*K)\impl$.
\end{proof}

The analogy between $\C^n$ and $(T^*K)\impl$ goes even further: if
$\Xi$ is integral, it generates a circle action, which turns out to be
locally free.  This will be proved in greater generality in the next
section.

\section{The stratification}\label{section;stratification}

We show now that the symplectic decomposition of the imploded
cross-section of a Hamiltonian $K$-manifold $M$ is a
\emph{stratification} in the sense that it is locally finite,
satisfies the frontier condition (i.e.\ $X_i\cap\bar X_j\ne\emptyset$
implies $X_i\subset\bar X_j$) and a certain regularity condition,
which is sometimes called local normal triviality.  This means that
locally at every point, in the direction transverse to the stratum,
$M\impl$ is a cone over a lower-dimensional stratified space, called
the \emph{link} of the point.  The link carries a locally free circle
action such that the quotient space, the \emph{symplectic link},
decomposes naturally into symplectic manifolds.  (In fact, the
symplectic link is the imploded cross-section of a singular symplectic
quotient.)
 
This is analogous to results about singular symplectic quotients
proved in \cite[\S5]{sjamaar-lerman;stratified}.  (There is however
one aspect in which imploded cross-sections are different from
symplectic quotients: they do not appear to have a naturally defined
algebra of functions equipped with a Poisson bracket.  Imploded
cross-sections are therefore strictly speaking not ``stratified
symplectic spaces'' in the sense of \cite{sjamaar-lerman;stratified}.)
The strategy of the proof is the same: write a local normal form for
an open neighbourhood of any point in $\cross$ and carry out all
computations in this model.  We shall do this in three steps.  The
final result is summarized in Theorem \ref{theorem;local} at the end
of this section.  Let $x\in M\impl$ and choose $m\in\cross$ such that
$x=\pi(m)$.

\emph{Step} 1.  Assume that $K$ is semisimple and $\Phi(m)=0$.  Then
the orbit $Km$ is isotropic.  The \emph{symplectic slice} at $m$ is
$$
V=T_m(Km)^\omega\!/T_m(Km).
$$
The natural linear action of the isotropy subgroup $H=K_m$ on $V$ is
symplectic and has a moment map given by \eqref{equation;quadratic}.
Hence we can form the quotient
\begin{equation}\label{equation;local-model}
F(K,H,V)=(T^*K\times V)\quot H,
\end{equation}
where we let $H$ act on $T^*K$ from the left.  It is a symplectic
vector bundle with fibre $V$ over the base $T^*K\quot H=T^*(K/H)$, and
it carries a Hamiltonian $K$-action induced by the right $K$-action on
$T^*K$.  Let $[k,\lambda,v]\in F(K,H,V)$ denote the point
corresponding to a point $(k,\lambda,v)\in T^*K\times
V=K\times\k^*\times V$ in the zero fibre of the $H$-moment map.  The
symplectic slice theorem of Marle and Guillemin-Sternberg says that
there exists a map from a $K$-stable open neighbourhood of $m$ in $M$
to a $K$-stable open neighbourhood of $[1,0,0\,]$ in $F(K,H,V)$ which
is an isomorphism of Hamiltonian $K$-manifolds and sends $m$ to
$[1,0,0\,]$.  This means that for the purpose of investigating
$M\impl$ near $x=\pi(m)$ we can replace $M$ by $F(K,H,V)$.  Corollary
\ref{corollary;reduction-implosion} implies
$$
F(K,H,V)\impl=\bigl((T^*K)\impl\times V\bigr)\bigquot H.
$$
Since $K$ is semisimple, the lowest-dimensional stratum of
$(T^*K)\impl$ consists of a single point.  Therefore the stratum of
$\pi([1,0,0\,])$ in $F(K,H,V)\impl$ is $V^H$, which shows once again
that the stratum of $x$ in $M\impl$ is a locally closed subspace and a
symplectic manifold.  Let $W=(V^H)^\omega$ be the skew complement of
$V^H$.  Then $W^H=\{0\}$ and $V=W\oplus V^H$ as a symplectic
$H$-module.  Thus we have a splitting
\begin{multline*}
F(K,H,V)\impl=\bigl((T^*K)\impl\times W\times V^H\bigr)\bigquot H\\
=\bigl(((T^*K)\impl\times W)\quot H\bigr)\times V^H
=F(K,H,W)\impl\times V^H,
\end{multline*}
which shows that, near $x$, $M\impl$ is symplectically the product of
a neighbourhood of $x$ inside its stratum and a neighbourhood of
$\pi([1,0,0\,])$ inside $F(K,H,W)\impl$.

The space $F(K,H,W)\impl=\bigl((T^*K)\impl\times W\bigr)\bigquot H$
carries information on the nature of the singularity at $x$.  Let
\begin{equation}\label{equation;cone-decompose}
F(K,H,W)\impl=\coprod_{j\in\ca J}F_j
\end{equation}
be its symplectic decomposition, let $*=F_{j_0}=\pi([1,0,0\,])$ be the
lowest stratum, called the \emph{vertex}, and put
$$
F^\times(K,H,W)\impl=\coprod_{j\in\ca J\setminus\{j_0\}}\!\!\!\!F_j.
$$
Define a continuous $\R^\times$-action $\hat{\ca A}$ on
$(T^*K)\impl\times W$ by
$$
\hat{\ca A}_t(\pi_{\ca R}(k,\lambda),w) =\bigl(\pi_{\ca
R}(k,t\lambda),\sqrt{t}\,w\bigr).
$$
This action commutes with the $H$-action, preserves the strata of
$(T^*K)\impl\times W$, and on each stratum is a Liouville action.
Moreover, the $H$-moment map is homogeneous of degree $1$ with respect
to $\hat{\ca A}$, so $\hat{\ca A}$ descends to an action $\ca A$ on
$F(K,H,W)\impl$ which preserves the stratification and on each stratum
is a Liouville action.  Let $\omega_j$ be the symplectic form,
$\zeta_j$ the Liouville vector field, and $\beta_j$ the one-form
$\iota(\zeta_j)\omega_j$ on $F_j$.

Now choose an $H$-invariant $\omega$-compatible complex structure on
the $H$-module $W$ and a circle subgroup of $T$ with infinitesimal
generator $\Xi\in\t$.  Consider the $S^1$-action on $(T^*K)\impl\times
W$ given by $\Xi$ acting from the right on $(T^*K)\impl$ and complex
scalar multiplication on $W$.  It is generated by the Hamiltonian
$\hat\varphi(\pi_{\ca R}(k,\lambda),w)
=\lambda(\Xi)-\frac1{2}\norm{w}^2$, it commutes with the $H$-action,
and therefore descends to a Hamiltonian $S^1$-action on
$F(K,H,W)\impl$.  The reduced Hamiltonian is
$\varphi(\pi([k,\lambda,w])) =\lambda(\Xi)-\frac1{2}\norm{w}^2$.  For
each $j\in\ca J$ let $\Xi_j$ be the vector field on $F_j$ induced by
$\Xi\in\t$.  The functions $\varphi_j=\varphi|_{F_j}$ satisfy
$\varphi_j =-\iota(\Xi_j)\beta_j$.  The subspace
$\lk(x)=\varphi\inv(-1)$ of $F(K,H,W)\impl$ is called the \emph{link}
of $x$.  The following result says that, stratum by stratum,
$F(K,H,W)\impl$ is the symplectification of $\lk(x)$ in the sense of
Arnold \cite[Appendix 4]{arnold;mathematical-methods}.  It is proved
in the same way as Proposition \ref{proposition;homogeneous}.

\begin{proposition}\label{proposition;implode-homogeneous}
Assume that $K$ is semisimple and that $\Xi$ is in the interior of the
cone spanned by the negative coroots.  Then
\begin{enumerate}
\item\label{item;implode-proper}
$\varphi$ is proper on $F(K,H,W)\impl$\upn, $\lk(x)$ is compact\upn,
$\varphi\le0$\upn, and $\varphi\inv(0)=*$\upn;
\item\label{item;implode-subcontact}
for all $j\ne j_0$\upn, $\lk(x)_j=\lk(x)\cap F_j$ is a smooth
manifold\upn, $\beta_j$ restricts to a contact form\upn, and $\Xi_j$
to the Reeb vector field on $\lk(x)_j$\upn;
\item\label{item;implode-section}
the map $f\colon\lk(x)\times\R^\times\to F^\times(K,H,W)\impl$ defined
by
$$
f\bigl(\pi([k,\lambda,w]),t\bigr)=\ca
A_t\bigl(\pi([k,\lambda,w])\bigr)
$$
is a stratum-preserving homeomorphism\upn, which on every stratum is a
diffeomorphism\upn, and satisfies
$f^*\omega_j=d(t(\beta_j|_{\lk(x)_j}))$\upn;
\item\label{item;implode-cone}
$f$ extends uniquely to a homeomorphism $C^\circ(\lk(x))\to
F(K,H,W)\impl$.
\end{enumerate}
\end{proposition}

Consequently, if the link is a (homology) sphere, then $M\impl$ is a
topological (homology) manifold at $x$.  The quotient
$$
\slk(x)=\lk(x)/S^1=F(K,H,W)\impl\quot[-1]S^1 =\bigl((T^*K)\impl\times
W\bigr)\bigquot[(0,-1)]H\times S^1
$$
is the \emph{symplectic link} of $x$.  For instance, the symplectic
link of the vertex in $(T^*K)\impl$ is $(T^*K)\impl\quot[-1]S^1$.  To
be definite, let us henceforth take $\Xi$ to be the sum of the simple
coroots.  (This choice is motivated by Proposition
\ref{proposition;slink} in the next section.)  The following result
says that, stratum by stratum, the link is the contactification of the
symplectic link.

\begin{proposition}\label{proposition;circle}
Assume that $K$ is semisimple.  Then
\begin{enumerate}
\item\label{item;locally-free}
the $S^1$-action on $F^\times(K,H,W)\impl$ is locally free.
\item\label{item;orbifold}
for each $j\in\ca J\setminus\{j_0\}$ the space
$\slk(x)_j=\lk(x)_j/S^1$ is a symplectic orbifold\upn, and
$\beta_j|_{\lk(x)_j}$ is a connection on the orbibundle
$\lk(x)_j\to\slk(x)_j$\upn, whose curvature is the reduced symplectic
form on $\slk(x)_j$.
\end{enumerate}
\end{proposition}

\begin{proof}
Let us denote the element $(\pi_{\ca R}(k,\lambda),w)$ of
$(T^*K)\impl\times W$ by $(\!(k,\lambda,w)\!)$ and its image in the
orbit space $\bigl((T^*K)\impl\times W\bigr)\big/S^1$ by $\powl
k,\lambda,w\powr$.  To prove \eqref{item;locally-free} we need to show
that for every $(\!(k,\lambda,w)\!)\ne(\!(1,0,0)\!)$ in the zero level
set of the $H$-moment map, the infinitesimal stabilizers $\lie
l_{(\!(k,\lambda,w)\!)}$ and $\lie l_{\powl k,\lambda,w\powr}$ are
equal.  There exists an infinitesimal character $\chi\in\lie l_{\powl
k,\lambda,w\powr}^*$ such that for all $\eta\in\lie l_{\powl
k,\lambda,w\powr}$
\begin{equation}\label{equation;character}
\exp\eta\cdot(\!(k,\lambda,w)\!)  =\bigl(\!\bigl(k\exp(-2\pi
i\chi(\eta)\Xi),\lambda,e^{2\pi i\chi(\eta)}w\bigr)\!\bigr).
\end{equation}
Let $\sigma$ be the face containing $\lambda$.  Then
\eqref{equation;character} boils down to
$$
\exp(k\inv\eta)=\exp(-2\pi i\chi(\eta)\Xi)\bmod[K_\sigma,K_\sigma]
\quad\text{and}\quad(\exp\eta)w=e^{2\pi i\chi(\eta)}w.
$$
Differentiating at $\eta=0$ yields
\begin{align}
\label{equation;character1}
k\inv\eta+\chi(\eta)\Xi&=[\xi_1,\xi_2]\quad\text{for certain $\xi_1$,
$\xi_2\in\k_\sigma$},\\
\label{equation;character2}
\eta w&=\chi(\eta)w.
\end{align}
Since $(\!(k,\lambda,w)\!)$ is in the zero fibre of the $H$-moment
map,
$$
0=\Phi_{\ca L}^\eta(k,\lambda)+\Phi_W^\eta(w)
=-\lambda(k\inv\eta)+\frac1{2}\omega_W(\eta w,w).
$$
Combined with \eqref{equation;character2} this gives
$$
\lambda(k\inv\eta) =\frac1{2}\omega_W(\eta w,w)
=\frac{\chi(\eta)}{2}\omega_W(w,w)=0.
$$
Applying $\lambda$ to both sides of \eqref{equation;character1} we
then obtain
$$
\chi(\eta)\lambda(\Xi)=\lambda([\xi_1,\xi_2])
=-(\ad^*(\xi_1)\lambda)(\xi_2)=0,
$$
because $\xi_1\in\k_\sigma$ and $\lambda\in\sigma$.  Since $\Xi<0$ on
$\t^*_+$, $\lambda(\Xi)<0$ and hence $\chi(\eta)=0$.  From
\eqref{equation;character} we conclude that $\eta\in\lie
l_{(\!(k,\lambda,w)\!)}$, and therefore $\lie l_{\powl
k,\lambda,w\powr}$ is contained in $\lie l_{(\!(k,\lambda,w)\!)}$.
The reverse inclusion is obvious.  This proves
\eqref{item;locally-free}.  \eqref{item;orbifold} follows immediately
from \eqref{item;locally-free}.
\end{proof}
 
Theorem \ref{theorem;decomposition} implies that the symplectic link
is a union of symplectic manifolds.  In fact, the manifold
decomposition of $\slk(x)$ is a refinement of the orbifold
decomposition $\slk(x) =\coprod\slk(x)_j$.

\emph{Step} 2.  Next consider the case where $K$ may have a
positive-dimensional centre and $\Phi(m)$ is contained in $\z^*$, the
fixed point set of the coadjoint action.  This case reduces
immediately to the previous one by replacing $K$ with its semisimple
part $[K,K]$.  This works because the $[K,K]$-moment map sends $m$ to
$0\in[\k,\k]^*$ and if $U$ is a $[K,K]$-stable open neighbourhood of
$m$, then $U\cap\cross$ is saturated under the equivalence relation
$\sim$, so that $\pi(U\cap\cross)$ is an open neighbourhood of
$\pi(m)$ in $M\impl$.

\emph{Step} 3.  Finally we reduce the general case to the case
$\Phi(m)\in\z^*$.  Let $\sigma\in\Sigma$ be the face containing
$\Phi(m)$.  Then $m\in M_\sigma$ and $\pi(m)\in
M_\sigma\quot{[K_\sigma,K_\sigma]}$.  The \emph{standard open
neighbourhood} of $M_\sigma\quot{[K_\sigma,K_\sigma]}$ in $M\impl$ is
the $T$-stable open set
$$
O_\sigma=\Phi\inv(\star\sigma)/{\sim} =\Phi\impl\inv(\star\sigma)
=\coprod_{\tau\ge\sigma}M_\tau\quot{[K_\tau,K_\tau]}.
$$
Let $R\subset\t^*$ be the root system of $(K,T)$ and $S$ the set of
roots which are simple relative to the chamber $\t^*_+$.  The root
system of $(K_\sigma,T)$ is then
$$
R_\sigma=\{\,\alpha\in R\mid\text{$\lambda(\alpha\spcheck)=0$ for all
$\lambda\in\sigma$}\,\},
$$ 
and its set of simple roots is $S_\sigma=R_\sigma\cap S$.  The
corresponding positive Weyl chamber is denoted by $\t^*_{+,\sigma}$.
Both $\z^*_\sigma$ and $\star\sigma$ are contained in
$\t^*_{+,\sigma}$.  Let
$(M_\sigma)\impl=\Phi_\sigma\inv(\t^*_{+,\sigma})/{\sim_\sigma}$ be
the imploded cross-section of $M_\sigma$ with respect to the
$K_\sigma$-action,
$\pi_\sigma\colon\Phi_\sigma\inv(\t^*_{+,\sigma})\to(M_\sigma)\impl$
the quotient map, and $(\Phi_\sigma)\impl$ the associated imploded
moment map.

\begin{lemma}\label{lemma;open}
\begin{enumerate}
\item\label{item;open}
$\star\sigma$ is open in $\t^*_{+,\sigma}$.
\item\label{item;standard-open}
$O_\sigma$ is isomorphic to the open subset
$(\Phi_\sigma)\impl\inv(\star\sigma)$ of $(M_\sigma)\impl$.
\end{enumerate}
\end{lemma}

\begin{proof}
\eqref{item;open} follows immediately from
\begin{align*}
\t^*_{+,\sigma}
&=\bigl\{\,\lambda\in\t^*:\text{$\lambda(\alpha\spcheck)\ge0$ for all
$\alpha\in S_\sigma$}\,\bigr\},\\
\star\sigma
&=\bigl\{\,\lambda\in\t^*:\text{$\lambda(\alpha\spcheck)\ge0$ for all
$\alpha\in S_\sigma$, $\lambda(\alpha\spcheck)>0$ for all $\alpha\in
S\setminus S_\sigma$}\,\bigr\}.
\end{align*}

From \eqref{item;open} it follows that
$(\Phi_\sigma)\impl\inv(\star\sigma)$ is open in $(M_\sigma)\impl$.
For $\lambda\in\star\sigma$, $K_\lambda\subset K_\sigma$, so for
$m_1$, $m_2\in\Phi\inv(\star\sigma)$, $m_1\sim m_2$ is equivalent to
$m_1\sim_\sigma m_2$.  Hence $(\Phi_\sigma)\impl\inv(\star\sigma)
=\Phi_\sigma\inv(\star\sigma)/{\sim_\sigma}
=\Phi\inv(\star\sigma)/{\sim}=O_\sigma$.
\end{proof}

To examine the structure of $M\impl$ near $\pi(m)$ we can therefore
resort to the space $(M_\sigma)\impl$ and the point $\pi_\sigma(m)$.
But here the argument of steps 1 and 2 applies, because
$\Phi(m)\in\z_\sigma^*$ is fixed under the coadjoint action of
$K_\sigma$.  We can summarize this discussion as follows.

\begin{theorem}\label{theorem;local}
Let $(M,\omega,\Phi)$ be a Hamiltonian $K$-manifold and let $x\in
M\impl$.  Choose $m\in\cross$ such that $x=\pi(m)$.  Let $\sigma$ be
the face of $\t^*_+$ containing $\Phi(m)$.  Let
$H=[K_\sigma,K_\sigma]_m$ be the stabilizer of $m$ and $V$ the
symplectic slice at $m$ for the $[K_\sigma,K_\sigma]$-action on
$M_\sigma$.  Put $W=(V^H)^\omega$.  Then $x$ has an open neighbourhood
which is isomorphic to a product $U_1\times U_2$\upn, where $U_1$ is a
neighbourhood of $x$ in its stratum and $U_2$ is a neighbourhood of
the vertex in $F([K_\sigma,K_\sigma],H,W)\impl$.  The space
$F([K_\sigma,K_\sigma],H,W)\impl$ is the stratified symplectification
of the link $\lk(x)$ in the sense of Proposition
{\rm\ref{proposition;implode-homogeneous}}\upn, and $\lk(x)$ is the
stratified contactification of the symplectic link $\slk(x)$ in the
sense of Proposition {\rm\ref{proposition;circle}}.  The symplectic
decomposition of $M\impl$ is locally finite and satisfies the frontier
condition.
\end{theorem}

\begin{proof}
Everything has been proved except the frontier condition.  Let
$X_{i_0}$ be the stratum of $x$ and suppose that $x$ is in the closure
of a stratum $X_i$.  Put $Y=\{\,y\in X_{i_0}\mid y\in\bar X_i\,\}$.
Then $x\in Y$ and $Y$ is closed in $X_{i_0}$.  In the local model
$F([K_\sigma,K_\sigma],H,V)\impl$ around $x$, the stratum $X_{i_0}$ is
of the form $F_{j_0}\times V^H=\{*\}\times V^H$, whereas $X_i$ is of
the form $F_j\times V^H$ for some $j\in\ca J$.  Here the notation is
as in \eqref{equation;cone-decompose}.  Every stratum $F_j$ in
$F^\times(K,H,W)\impl$ is stable under the $\R^\times$-action and
therefore has the vertex $*$ as a limit point.  It follows that
$X_{i_0}\cap U\subset\bar X_i\cap U$, where $U$ is an appropriate open
neighbourhood of $x$ in $M\impl$.  Hence $Y$ is open, and therefore
$Y=X_{i_0}$.  We have shown that $X_{i_0}\subset\bar X_i$.
\end{proof}


In the next section we shall prove that the link and the symplectic
link of every point in $M\impl$ are connected.

\section{K\"ahler structures}\label{section;kaehler}

In this section we show that the strata of the universal imploded
cross-section fit together in an unexpectedly nice way.  It turns out
that if $K$ is semisimple and simply connected, $(T^*K)\impl$ can be
embedded into a unitary $K$-module $E$ in such a manner that the
symplectic form on each stratum is the pullback of the flat K\"ahler
form on $E$.  The image of the embedding is a closed $K$-stable affine
subvariety with coordinate ring $\C[G]^N$, where $G=K^\C$ is the
complexification of $K$ and $N$ is a maximal unipotent subgroup of
$G$.  This happy state of affairs permits us to calculate
$(T^*K)\impl$ for groups other than $\SU(2)$ and to prove that the
imploded cross-section of \emph{every} K\"ahler Hamiltonian
$K$-manifold is a K\"ahler space.

Assume first that $K$ is an arbitrary compact connected Lie group.
Let $\Lambda=\ker({\exp}|_\t)$ be the exponential lattice in $\t$ and
$\Lambda^{\!*}=\Hom_\Z(\Lambda,\Z)$ the weight lattice in $\t^*$.
Then $\Lambda^{\!*}_+=\Lambda^{\!*}\cap\t^*_+$ is the monoid of
dominant weights.  For any dominant weight $\lambda$ let $V_\lambda$
be the $K$-module with highest weight $\lambda$.  Select a minimal set
of generators $\Pi$ of $\Lambda^{\!*}_+$, put
\begin{equation}\label{equation;bigmodule}
E=\bigoplus_{\varpi\in\Pi}V_\varpi,
\end{equation}
and fix a highest-weight vector $v_\varpi$ in each $V_\varpi$.  To
each face $\sigma\in\Sigma$ is associated a parabolic subgroup
$P_\sigma$ of $G$ with Lie algebra
$$
\lie p_\sigma=\t^\C\oplus\bigoplus_{\alpha\in
R_\sigma}\g_\alpha\oplus\bigoplus_{R_+\setminus
R_{+,\sigma}}\!\!\g_\alpha,
$$
where $\g_\alpha$ denotes the root space for $\alpha$.

\begin{lemma}\label{lemma;stabilizer}
For each face $\sigma$ let
$v_\sigma=\sum_{\varpi\in\bar\sigma}v_\varpi$.  Then the stabilizer of
$v_\sigma$ for the $G$-action is equal to $[P_\sigma,P_\sigma]$.
\end{lemma}

\begin{proof}
Let $[v_\sigma]$ denote the equivalence class of $v_\sigma$ in the
product of projective spaces
$\prod_{\varpi\in\bar\sigma}\P(V_\varpi)$.  Then
$G_{[v_\sigma]}=P_\sigma$ and
$$
G_{v_\sigma}=\{\,g\in P_\sigma\mid\text{$\lambda(g)=1$ for all
$\lambda\in\Lambda^{\!*}\cap\bar\sigma$}\,\}.
$$
(Here we identify $\lambda$ with the character of $P_\sigma$ that it
exponentiates to.)  Let $G_\sigma=(K_\sigma)^\C$ and let $U_\sigma$ be
the unipotent radical of $P_\sigma$.  The corresponding Lie algebras
are
\begin{equation}\label{equation;lie}
\g_\sigma=\t^\C\oplus\bigoplus_{\alpha\in R_\sigma}\g_\alpha
\qquad\text{and}\qquad\lie u_\sigma=\bigoplus_{\alpha\in R_+\setminus
R_{+,\sigma}}\g_\alpha,
\end{equation}
and we have Levi decompositions
\begin{equation}\label{equation;levi}
P_\sigma=G_\sigma U_\sigma\qquad\text{and}\qquad
[P_\sigma,P_\sigma]=[G_\sigma,G_\sigma]\,U_\sigma.
\end{equation}
Since characters of $P_\sigma$ vanish on $U_\sigma$, we have
$G_{v_\sigma} =Q_\sigma U_\sigma$ with $Q_\sigma=G_{v_\sigma}\cap
G_\sigma$.  Thus both $[G_\sigma,G_\sigma]$ and $Q_\sigma$ are
semisimple groups with root system $R_\sigma$.  To finish the proof it
suffices to show that $Q_\sigma$ is connected.  The elements of
$\Lambda^{\!*}\cap\bar\sigma$ being $G_\sigma$-invariant, $Q_\sigma$
is stable under conjugation by the connected group $G_\sigma$.
Therefore we need only show that the intersection of $Q_\sigma$ with
the maximal torus $T$ is connected.  Observe that
$$
Q_\sigma\cap T=\{\,t\in T\mid\text{$\lambda(t)=1$ for all
$\lambda\in\Lambda^{\!*}\cap\bar\sigma$}\,\}.
$$
If $C$ is any top-dimensional polyhedral cone in a vector space $V$
and $\Gamma$ is a lattice in $V$, then the set $\Gamma\cap C$ contains
a $\Z$-basis for $\Gamma$.  Therefore $\Lambda^{\!*}\cap\bar\sigma$
contains a $\Z$-basis for the lattice $\Lambda^{\!*}\cap\z_\sigma^*$.
This basis can be extended to a basis of $\Lambda^{\!*}$, and this
implies that $Q_\sigma\cap T$ is connected.
\end{proof}

Let $N$ be the maximal unipotent subgroup of $G$ with Lie algebra
$\lie n=\bigoplus_{\alpha\in R_+}\g_\alpha$ and let $A$ be the (real)
subgroup with Lie algebra $\lie a=i\t$.  Then $T^\C=T\!A$ and
$B=T\!AN$ is a Borel subgroup of $G$.  The Borel-Weil Theorem implies
that the ring of $N$-invariants is a multi-graded direct sum,
$$
\C[G]^N=\bigoplus_{\lambda\in\Lambda^{\!*}_+}V_\lambda.
$$
In particular it is generated by the finite-dimensional subspace
$E^*$.  It follows that there exists a closed $G$-equivariant
algebraic embedding $G_{\!N}\hookrightarrow E$.  (Following Kraft
\cite[Kapitel III]{kraft;geometrische}, for any affine $G$-variety $X$
we denote by $X_N$ the affine variety with coordinate ring $\C[X]^N$.)
By results of Vinberg and Popov
\cite[\S3]{vinberg-popov;class-quasihomogeneous}, $G_{\!N}$ consists
of finitely many $G$-orbits, which are labelled by the faces of the
cone $\t^*_+$.  In addition, $G_{\!N}$ contains the subspace $E^N$ and
$G_{\!N}=GE^N$.  The stabilizer of the orbit corresponding to the face
$\sigma$ is the group of all $g\in P_\sigma$ such that $\lambda(g)=1$
for all $\lambda\in\Lambda^{\!*}\cap\bar\sigma$, which is equal to
$[P_\sigma,P_\sigma]$ by Lemma \ref{lemma;stabilizer}.  In particular,
the open orbit is of type $G/N$.  Moreover, the embedding
$G_{\!N}\hookrightarrow E$ is uniquely determined by sending $1\bmod
N$ to the sum of the highest-weight vectors
$\sum_{\varpi\in\Pi}v_\varpi$.  We shall identify $G_{\!N}$ with its
image in $E$.  We turn $E$ into a $K\times T$-module by letting $T$
act on $V_\varpi$ with weight $-\varpi$.  Observe that a different
choice of highest-weight vectors leads to a new embedding
$G_{\!N}\hookrightarrow E$ which differs from the old by
multiplication by an element of the complex torus $T\!A$.

Let $\langle\cdot,\cdot\rangle$ be the unique $K$-invariant Hermitian
inner product on $E$ satisfying $\norm{v_p}=1$ for all $p$.  We regard
$E$ as a flat K\"ahler manifold with the K\"ahler form
$\omega_E=-\Im\langle\cdot,\cdot\rangle$.  It is convenient to write
$\omega_E=d\beta_E$ with
\begin{equation}\label{equation;connection}
(\beta_E)_v(w)=-\frac1{2}\Im\langle v,w\rangle
\end{equation}
for $v$, $w\in E$.

Now assume that $K$ is semisimple and simply connected.  The set $\Pi$
is then uniquely determined: it is the set of fundamental weights
$\{\varpi_1,\varpi_2,\dots,\varpi_r\}$, which form a $\Z$-basis of the
weight lattice.  Let $\alpha_1$, $\alpha_2,\dots$, $\alpha_r$ be the
corresponding simple roots.  Put $V_p=V_{\varpi_p}$ and
$v_p=v_{\varpi_p}$.  Since $\lambda(\alpha\spcheck)\ge0$ if
$\lambda\in\t^*_+$,
\begin{equation}\label{equation;embedding}
s(\lambda)=\frac1{\sqrt{\pi}}
\sum_{p=1}^r\sqrt{\lambda(\alpha\spcheck_p)}\,v_p
\end{equation}
defines a continuous map from $\t^*_+$ into $E^N$.

\begin{remark}\label{remark;toric}
The subspace $E^N$ inherits a symplectic form from $E$ and according
to \eqref{equation;quadratic} the moment map for the $T$-action on
$E^N$ is given by
$\Phi_{E^N}(\sum_pc_pv_p)=-\pi\sum_p\abs{c_p}^2\varpi_p$.  Observe
that, $K$ being semisimple and simply connected, the $T$-action on
$E^N$ is equivalent to the $(\C^\times)^r$-action on $\C^r$ and so is
effective and multiplicity-free.  The map $\Phi_{E^N}$ separates the
$T$-orbits and its image is the opposite chamber $-\t^*_+$.  Therefore
$E^N$ is nothing but the symplectic toric manifold (multiplicity-free
$T$-manifold or Delzant space) associated with the polyhedron
$-\t^*_+$.  Finally note that $\Phi_{E^N}(s(\lambda))=-\lambda$, that
is to say, $s$ is a section of $-\Phi_{E^N}$.
\end{remark}

The map $s$ extends uniquely to a $K\times T$-equivariant map $\ca
F\colon K\times\t^*_+\to E$.

\begin{proposition}\label{proposition;embedding}
Assume that $K$ is semisimple and simply connected.
\begin{enumerate}
\item\label{item;embedding}
$\ca F$ induces a map $f\colon(T^*K)\impl\to E$ which is continuous
and closed \upn(for the classical topology on $E$\upn)\upn, and
injective.
\item\label{item;symplectic}
The restriction of $f$ to each stratum is a smooth symplectomorphism.
\item\label{item;algebraic}
The image of $(T^*K)\impl$ under $f$ is identical to $G_{\!N}$.  Thus
$f\colon(T^*K)\impl\to G_{\!N}$ is an isomorphism of Hamiltonian
$K$-spaces in the sense of Section {\rm\ref{section;abelian}}.
\end{enumerate}
\end{proposition}

\begin{proof}
It is plain from \eqref{equation;embedding} that $\ca F$ is continuous
and closed.  Furthermore, by Lemma \ref{lemma;stabilizer} the
stabilizer of $\ca F(1,\lambda)=s(\lambda)$ for the $K$-action is
equal to $K\cap[P_\sigma,P_\sigma]=[K_\sigma,K_\sigma]$, where
$\sigma$ is the face containing $\lambda$.  This implies
\eqref{item;embedding}.

It is clear from \eqref{equation;embedding} that $\ca F$ is smooth on
$K\times\sigma$ for every $\sigma$, and therefore $f$ restricted to
$K/[K_\sigma,K_\sigma]\times\sigma$ is a smooth embedding.  We check
that it preserves the symplectic form by showing that
$f^*\beta_E=\beta_\sigma$, where $\beta_\sigma$ is the one-form on
$K/[K_\sigma,K_\sigma]\times\sigma$ considered in Lemma
\ref{lemma;exact}.  Because $\ca F$ is $K$-equivariant we need only
show this at the points $(\bar1,\lambda)$, where $\bar k$ denotes the
coset $k\,[K_\sigma,K_\sigma]$.  By Lemma \ref{lemma;exact},
\begin{equation}\label{equation;oneform}
(\beta_\sigma)_{(\bar1,\lambda)}(\bar\xi,\mu) =\lambda(\xi)
\end{equation}
for all $\xi\in\k$ and $\mu\in\z_\sigma$.  On the other hand,
$$
(f^*\beta_E)_{(\bar1,\lambda)}(\bar\xi,\mu)
=(\beta_E)_{f(\bar1,\lambda)}(f_*(\bar\xi,\mu))
=(\beta_E)_{s(\lambda)}(\ca F_*(\xi,\mu)).
$$
Here
$$
\ca F_*(\xi,\mu)=\frac{d}{dt}\exp(t\xi)s(\lambda+t\mu)\biggm|_{t=0}
=\xi_E(s(\lambda))+\sum_{p=1}^r\frac{\mu(\alpha\spcheck_p)}%
{2\sqrt{\pi\,\lambda(\alpha\spcheck_p)}}\,v_p.
$$
Together with \eqref{equation;connection} this yields
\begin{align*}
(f^*\beta_E)_{(\bar1,\lambda)}(\bar\xi,\mu)
&=\frac1{2}\Im\biggl\langle\xi_E(s(\lambda))
+\sum_p\frac{\mu(\alpha\spcheck_p)}%
{2\sqrt{\pi\,\lambda(\alpha\spcheck_p)}}\,
v_p,s(\lambda)\biggr\rangle\\
&=\frac1{2\pi}
\Im\biggl\langle\sum_p\sqrt{\lambda(\alpha\spcheck_p)}\,\xi_E(v_p),
s(\lambda)\biggr\rangle
+\frac1{2\pi}\Im\sum_p\frac{\mu(\alpha\spcheck_p)}{2}\\
&=\sum_p\lambda(\alpha\spcheck_p)\varpi_p(\xi)=\lambda(\xi),
\end{align*}
where we have used $\langle\xi_E(v_p),v_q\rangle=2\pi
i\varpi_p(\xi)\delta_{pq}$.  Comparing with \eqref{equation;oneform}
we conclude that $f^*\beta_E=\beta_\sigma$.  This proves
\eqref{item;symplectic}.

\eqref{item;algebraic} is a consequence of the Iwasawa decomposition.
Put $\lie a=i\t$, $A=\exp\lie a$ and, for each face $\sigma$, $\lie
n_\sigma=\bigoplus_{\alpha\in R_{+,\sigma}}\g_\alpha$ and
$N_\sigma=\exp\lie n_\sigma$.  Then $G=K\!AN$ and $G_\sigma=K_\sigma
AN_\sigma$.  Recall also that $P_\sigma=G_\sigma U_\sigma$.  Here
$G_\sigma$ and $U_\sigma$ are as in \eqref{equation;lie}.  Let
$$
\lie a_\sigma=\lie a\cap[\g_\sigma,\g_\sigma]=\bigoplus_{\alpha\in
S_\sigma}i\R\alpha\spcheck\quad\text{and}\quad\lie a_\sigma^\perp
=i\z\oplus\bigoplus_{\alpha\in S\setminus
S_\sigma}\!\!i\R\alpha\spcheck,
$$
so that $\lie a=\lie a_\sigma\oplus\lie a_\sigma^\perp$.  Writing
$A_\sigma=\exp\lie a_\sigma$ and $A_\sigma^\perp =\exp\lie
a_\sigma^\perp$ we find $A=A_\sigma\times A_\sigma^\perp$ and, using
\eqref{equation;levi},
$$
[P_\sigma,P_\sigma]=[G_\sigma,G_\sigma]\,U_\sigma
=[K_\sigma,K_\sigma]A_\sigma N_\sigma
U_\sigma=[K_\sigma,K_\sigma]A_\sigma N.
$$
Hence
$$
G/[P_\sigma,P_\sigma]
=K\!AN\big/\bigl([K_\sigma,K_\sigma]A_\sigma^\perp N\bigr)\cong
K/[K_\sigma,K_\sigma]\times A_\sigma^\perp
$$
as smooth $K$-manifolds.  To finish the proof it suffices to show that
$s(\sigma)$ is equal to the $A_\sigma^\perp$-orbit through $v_\sigma$
for all $\sigma$.  For $\lambda\in\t^*_+$ put
$$
\psi(\lambda)=\frac1{4\pi i}\sum_p
\bigl(\log\lambda(\alpha\spcheck_p)-\log\pi\bigr)
\alpha\spcheck_p\in\lie a,
$$
where the sum is over all $p$ such that
$\lambda(\alpha\spcheck_p)\ne0$.  For each face $\sigma$, $\psi$
defines a diffeomorphism from $\sigma$ to $\lie a_\sigma^\perp$, and
therefore ${\exp}\circ\psi\colon\sigma\to A_\sigma^\perp$ is also a
diffeomorphism.  Moreover, for $\lambda\in\sigma$
\begin{align*}
\exp\psi(\lambda)\cdot v_\sigma &=\sum_{\substack{p\\ \alpha_p\in
S\setminus S_\sigma}}\exp\biggl(\frac1{4\pi i}\sum_{\substack{q\\
\alpha_q\in S\setminus S_\sigma}}
\bigl(\log\lambda(\alpha\spcheck_q)-\log\pi\bigr)
\alpha\spcheck_q\biggr)\cdot v_p\\
&=\sum_p\exp\biggl(\frac1{2}\sum_q
\bigl(\log\lambda(\alpha\spcheck_q)-\log\pi\bigr)
\varpi_p(\alpha\spcheck_q)\biggr)\cdot v_p\\
&=\sum_p\exp\frac1{2}
\bigl(\log\lambda(\alpha\spcheck_p)-\log\pi\bigr)\cdot v_p =\ca
F(\lambda).
\end{align*}
Hence $s(\sigma)=A_\sigma^\perp v_\sigma$.
\end{proof}

In a similar fashion the symplectic link of the vertex in
$(T^*K)\impl$ can be identified with a \emph{projective} variety.
Observe that the subvariety $G_{\!N}$ of $E$ is conical, i.e.\
preserved by the standard $\C^\times$-action on $E$.  The easiest way
to see this is to consider the one-parameter subgroup of $T$ generated
by $\Xi=-\sum_{p=1}^r\alpha\spcheck_p\in\t$.  As $\varpi_p(\Xi)=-1$
for all $p$ and $T$ acts with weight $-\varpi_p$ on $V_p$, $\Xi$
generates the standard $S^1$-action on $E$.  Since $G_{\!N}$ is affine
and is preserved under the action of $T$, it is a cone.  Let us denote
the subvariety $(G_{\!N}\setminus\{0\})/\C^\times$ of $\P(V)$ by
$\P(G_{\!N})$.  As before, $*$ denotes the vertex in $(T^*K)\impl$,
and $\slk(*) =(T^*K)\impl\quot[-1]S^1$ its symplectic link.  The
following result is now immediate from Proposition
\ref{proposition;embedding} and Example \ref{example;projective}.

\begin{proposition}\label{proposition;slink}
Assume that $K$ is semisimple and simply connected.  The isomorphism
$f\colon(T^*K)\impl\to G_{\!N}$ induces an isomorphism of Hamiltonian
$K$-spaces $\slk(*)\cong\P(G_{\!N})$.
\end{proposition}

Now let $K$ be a torus.  Then $\Lambda^{\!*}_+=\Lambda^{\!*}$ and
$(T^*K)\impl=T^*K$.  Let us take
$\Pi=\{\pm\kappa_1,\pm\kappa_2,\dots,\pm\kappa_s\}$, where
$\{\kappa_1,\kappa_2,\dots,\kappa_s\}$ is a $\Z$-basis of the lattice
$\Lambda^{\!*}$.  Let $\eta_1$, $\eta_2,\dots$, $\eta_s$ be the dual
basis of $\Lambda$, $V_p=V_{\kappa_p}$, and $V_{-p}=V_{-\kappa_p}$.
For $\lambda\in\k^*$ set
$$
s(\lambda)=\frac1{\sqrt{2\pi}}
\sum_{p=1}^s\biggl(\chi(\lambda(\eta_p))\,v_p
+\frac1{\chi(\lambda(\eta_p))}\,v_{-p}\biggr)
$$
with $\chi(t)=\sqrt{t+\sqrt{t^2+1}}$.  Then $s$ extends uniquely to an
equivariant map $\ca F$ from $T^*K$ into $E$ and it is straightforward
to check that this is a symplectic embedding.

The main points of this discussion can be restated as follows.

\begin{theorem}\label{theorem;affine}
Assume that $K$ is the product of a torus and a semisimple simply
connected group.  There exists a $K\times T$-equivariant embedding $f$
of $(T^*K)\impl$ into the unitary $K\times T$-module $E$ whose image
is the Zariski-closed affine subvariety $G_{\!N}$.  Hence the action
of $K\times T$ on $(T^*K)\impl$ extends to an action of the
complexified group $K^\C\times T^\C=G\times T\!A$.  The strata of
$(T^*K)\impl$ coincide with the orbits of $G$\upn:
$$
f\bigl((K\times\eu S_\sigma)\quot{[K_\sigma,K_\sigma]}\bigr)
=G/[P_\sigma,P_\sigma]
$$
for all faces $\sigma$.  The symplectic form on each stratum is the
restriction of the flat K\"ahler form on $E$.
\end{theorem}

\begin{example}\label{example;su2continued}
Let $K=\SU(2)$.  Write an arbitrary element of $G=\SL(2,\C)$ as
$g=(x_{ij})$.  Let $N$ be the subgroup consisting of upper triangular
unipotent matrices.  The $N$-invariants of degree $1$ are the entries
in the first column, $x_{11}$ and $x_{21}$.  These two elements freely
generate $\C[G]^N$.  Therefore $G_{\!N}$ is the affine plane $\C^2$,
which confirms the computation in Example \ref{example;su2}.
\end{example}

For general compact connected $K$ there is a similar embedding of
$(T^*K)\impl$ into $E$, but we have not been able to find one that is
symplectic with respect to a natural K\"ahler structure on $E$.
Instead we proceed as follows.  Consider the universal cover
$[K,K]\sptilde$ of $[K,K]$ and the group $\tilde
K=Z\times[K,K]\sptilde$.  Then $K=\tilde K/\Upsilon$, where $\Upsilon$
is a finite central subgroup of $\tilde K$.  Let $\tilde G$ be the
complexification of $\tilde K$ and $\tilde N$ the preimage of $N$ in
$\tilde G$.  Then $(T^*K)\impl\cong(T^*\tilde K)\impl/\Upsilon$ (see
Example \ref{example;semisimple}) and likewise $G_{\!N}\cong\tilde
G_{\!\tilde N}/\Upsilon$.  It follows that $f$ descends to a
homeomorphism $(T^*K)\impl\to G_{\!N}$.  We use this map to identify
$(T^*K)\impl$ with $G_{\!N}$, thus defining a structure of an affine
variety on $(T^*K)\impl$ and K\"ahler structures on the orbits of
$G_{\!N}$.

By virtue of this result we can bring the machinery of algebraic
geometry to bear on the universal imploded cross-section.  For
instance, it now makes sense to talk about algebraic subvarieties of
$(T^*K)\impl$.  Each stratum, being an orbit of $G$, is a quasi-affine
subvariety and its closure in the classical topology is the same as
its Zariski closure.  The following is an algebraic slice theorem for
$G_{\!N}$, valid for arbitrary reductive $G$.

\begin{lemma}\label{lemma;slice}
For every face $\sigma$ the point $v_\sigma$ has a $G$-stable
Zariski-open neighbourhood in $G_{\!N}$ which is equivariantly
isomorphic to $G\times^{[P_\sigma,P_\sigma]}[P_\sigma,P_\sigma]_{N}$.
\end{lemma}

\begin{proof}
Let $E_\sigma=\bigoplus_{\varpi\in\bar\sigma}V_\varpi$ and let
$\pr\colon E\to E_\sigma$ be the orthogonal projection.  Then $\pr$ is
$G$-equivariant and $\pr(G_{\!N})$ is the closure of $Gv_\sigma$.  For
any face $\tau$,
$\pr(v_\tau)=\sum_{\varpi\in\bar\sigma\cap\bar\tau}v_\varpi
=v_{\sigma\wedge\tau}$, where $\sigma\wedge\tau\in\Sigma$ is the
interior of $\bar\sigma\cap\bar\tau$.  Hence
\begin{equation}\label{equation;greater}
\pr(v_\tau)=v_\sigma\iff\tau\ge\sigma.
\end{equation}
Consider the subsets of $E$ given by
$$
X_\sigma=\coprod_{\tau\ge\sigma}[P_\sigma,P_\sigma]\,v_\tau
\qquad\text{and}\qquad
O_\sigma=GX_\sigma=\coprod_{\tau\ge\sigma}Gv_\tau.
$$
Then $X_\sigma$ is $[P_\sigma,P_\sigma]$-stable and
\eqref{equation;greater} implies that $X_\sigma =\pr\inv(v_\sigma)\cap
G_{\!N}$.  Hence $X_\sigma$ is Zariski-closed.  Similarly, $O_\sigma$
is equal to $\pr\inv(Gv_\sigma)\cap G_{\!N}$, and it is a $G$-stable
Zariski-open neighbourhood of $v_\sigma$.  If $x$ and $gx$ are in
$X_\sigma$, then $gv_\sigma=g\pr(x)=\pr(gx)=v_\sigma$, so that
$g\in[P_\sigma,P_\sigma]$ by Lemma \ref{lemma;stabilizer}.  It follows
that the multiplication map $G\times X_\sigma\to O_\sigma$ induces a
$G$-equivariant isomorphism $G\times^{[P_\sigma,P_\sigma]}X_\sigma\to
O_\sigma$.  The affine $[P_\sigma,P_\sigma]$-variety $X_\sigma$ is the
union of all orbits $[P_\sigma,P_\sigma]/[P_\tau,P_\tau]$ with
$\tau\ge\sigma$.  The groups $P_\tau$ are exactly the parabolic
subgroups of $P_\sigma$ that contain $B$, and therefore it follows
from the corollary to \cite[Theorem
6]{vinberg-popov;class-quasihomogeneous} that
$O_\sigma\cong[P_\sigma,P_\sigma]_{N}$.
\end{proof}

As an application we determine the smooth (nonsingular) locus of
$G_{\!N}$.  Since $G_{\!N}$ is smooth at $x$ if and only if it is
smooth at all points in $Gx$, it suffices to consider $x=v_\sigma$.

\begin{proposition}\label{proposition;locus}
Let $\sigma$ be any face of $\t^*_+$.
\begin{enumerate}
\item\label{item;smooth}
$G_{\!N}$ is smooth at $v_\sigma$ if and only if
$[G_\sigma,G_\sigma]\cong\SL(2,\C)^k$ for some $k$.  The slice
$[P_\sigma,P_\sigma]_N$ is then $\SL(2,\C)^k$-equivariantly isomorphic
to the standard $\SL(2,\C)^k$-representation on $(\C^2)^k$.
\item\label{item;finite-quotient}
$G_{\!N}$ has an orbifold singularity at $v_\sigma$ if and only if
$[G_\sigma,G_\sigma]\cong\SL(2,\C)^k/\Upsilon$ for some $k$ and some
central subgroup $\Upsilon$ of $\SL(2,\C)^k$.  Then
$[P_\sigma,P_\sigma]_N\cong(\C^2)^k/\Upsilon$ as a
$\SL(2,\C)^k/\Upsilon$-variety.
\end{enumerate}
\end{proposition}

\begin{proof}
It follows from Lemma \ref{lemma;slice} that $G_{\!N}$ is smooth at
$v_\sigma$ if and only if
$[P_\sigma,P_\sigma]_N\cong[G_\sigma,G_\sigma]_{N_\sigma}$ is smooth
at the vertex.  Pauer \cite[Hilfssatz 2.4]{pauer;glatte} proved that
this is the case if and only if $[G_\sigma,G_\sigma]\cong\SL(2,\C)^k$
and $[G_\sigma,G_\sigma]_{N_\sigma}\cong(\C^2)^k$.  This proves
\eqref{item;smooth}.

Likewise, $G_{\!N}$ has an orbifold singularity at $v_\sigma$ if and
only if $[G_\sigma,G_\sigma]_{N_\sigma}$ has an orbifold singularity
at the vertex.  This is the case precisely when $\ca G_{\ca N}$ has an
orbifold singularity at the vertex, where $\ca G$ is the universal
cover of $[G_\sigma,G_\sigma]$ and $\ca N$ is the preimage of
$N_\sigma$ in $\ca G$.  Since $\ca G$ is simply connected and $\ca N$
is connected, $\ca G/\ca N$ is simply connected.  The complement of
$\ca G/\ca N$ inside $(\ca G_{\ca N})\reg$ has complex codimension at
least $2$, and therefore $(\ca G_{\ca N})\reg$ is also simply
connected.  Moreover, $(\ca G_{\ca N})\reg$ is a union of $\ca
G$-orbits and therefore homogeneous, i.e.\ stable under dilations on
$\ca G_{\ca N}$.  We conclude that the vertex has a basis of open
neighbourhoods $O$ such that the nonsingular part $O\reg$ is simply
connected.  A result of Prill \cite[Theorem 1]{prill;local} implies
that if a complex analytic space $X$ has a finite-quotient singularity
at $x$ and if $x$ has a basis of neighbourhoods $O$ such that $O\reg$
is simply connected, then $X$ is smooth at $x$.  We conclude that $\ca
G_{\ca N}$ has an orbifold singularity at the vertex if and only if it
is smooth.  Therefore \eqref{item;finite-quotient} follows from
\eqref{item;smooth}.
\end{proof}

In particular $G_{\!N}$ is smooth if and only if $[G,G]$ is a product
of copies of $\SL(2,\C)$; it is an orbifold if and only if the
universal cover of $[G,G]$ is a product of copies of $\SL(2,\C)$.

\begin{example}\label{example;su3}
Let $K=\SU(3)$.  Let $(x_{ij})\in G=\SL(3,\C)$ and let $N$ be the
subgroup consisting of upper triangular unipotent matrices.  The
$N$-invariants of degree $1$ are the entries in the first column,
$w_1=x_{11}$, $w_2=x_{21}$, and $w_3=x_{31}$.  The invariants of
degree $2$ are the minors extracted from the first two columns,
$$
z_1=\begin{vmatrix}x_{21}&x_{22}\\x_{31}&x_{32}\end{vmatrix},\quad
z_2=\begin{vmatrix}x_{11}&x_{12}\\x_{31}&x_{32}\end{vmatrix},\quad
z_3=\begin{vmatrix}x_{11}&x_{12}\\x_{21}&x_{22}\end{vmatrix}.
$$
These six elements generate $\C[G]^N$ and the only relation is
$\sum_kw_kz_k=0$.  Thus $(T^*K)\impl=G_{\!N}$ is the quadric in
$\C^6=\C^3\times\C^3$ given by this equation.  The four strata
corresponding to the faces of $\t^*_+$ are $\{0\}$,
$(\C^3\setminus\{0\})\times\{0\}$, $\{0\}\times(\C^3\setminus\{0\})$
and the open stratum.  Only the vertex $\{0\}$ is singular.
\end{example}

As an application of the foregoing results let us show that the
imploded cross-section of an affine Hamiltonian $K$-space $X$ (as
defined in Example \ref{example;affine}) is an affine variety.  The
following result says that implosion is the symplectic analogue of
taking the quotient of a variety by the maximal unipotent subgroup of
a reductive group.  Recall that $X_N$ denotes the affine variety with
coordinate ring $\C[X]^N$.

\begin{theorem}\label{theorem;unipotent-algebraic}
The imploded cross-section of an affine Hamiltonian $K$-space $X$ is
$T$-equivariantly homeomorphic to the affine variety $X_N$.  Under the
homeomorphism the strata of $X\impl$ correspond to algebraic
subvarieties of $X_N$.
\end{theorem}

\begin{proof}
Embed $X$ into a finite-dimensional unitary $K$-module $V$ as in
Example \ref{example;affine}.  The $K$-action on $X$ extends uniquely
to an algebraic $G$-action, which preserves the stratification of $X$.
To examine $X\impl$ let us assume that $K$ is semisimple and simply
connected.  (This is justified by Lemma \ref{lemma;semisimple}.)  Then
$X\times G_{\!N}$ is a closed $K$-stable affine subvariety of $V\times
E$, and Theorems \ref{theorem;universal} and \ref{theorem;affine} give
us isomorphisms
\begin{equation}\label{equation;git}
X\impl\cong(X\times(T^*K)\impl)\quot K\cong(X\times G_{\!N})\quot K.
\end{equation}
A well-known result of Kempf and Ness \cite{kempf-ness;length} (see
also Schwarz \cite{schwarz;topology-algebraic}) says that the
symplectic quotient on the right is homeomorphic to the
invariant-theoretic quotient of $X\times G_{\!N}$ by $G$, i.e.\ the
affine variety (associated to the scheme) $\Spec\C[X\times
G_{\!N}]^G$.  The homeomorphism is induced by the inclusion of the
zero fibre of the moment map for the $K$-action on $X\times G_{\!N}$
and is therefore $T$-equivariant.  According to Kraft
\cite[III.3.2]{kraft;geometrische} the ring $\C[X\times G_{\!N}]^G$ is
isomorphic to $\C[X]^N$.
\end{proof}


An example of an affine Hamiltonian $K$-space is the local model space
defined in \eqref{equation;local-model}.  Thus we see that, for every
Hamiltonian $K$-manifold $M$, every point $x$ in $M\impl$ has an open
neighbourhood which is homeomorphic to an open subset (in the
classical topology) of an affine variety of the form $X_N$, where $X$
is the local model at $x$.  Since $X$ is smooth, the quotient $X_N$ is
normal (see e.g.\ Satz 2 in Kraft \cite[III.3.4]{kraft;geometrische}),
from which it follows that the link of $x$ is connected.

\begin{corollary}
Let $M$ be an arbitrary Hamiltonian $K$-manifold.  Then the link and
the symplectic link of every point in $M\impl$ are connected.
\end{corollary}

An assertion comparable to Theorem \ref{theorem;unipotent-algebraic}
can be made in the analytic category.  Let $(M,\omega,\Phi)$ be a
Hamiltonian $K$-manifold equipped with a $K$-invariant complex
structure $J$.  We assume that $J$ is compatible with $\omega$, so
that $M$ is a K\"ahler manifold, and that the $K$-action extends to a
holomorphic $G$-action.  We wish to show that $M\impl$ is a K\"ahler
space and in particular that its strata are K\"ahler manifolds.
However, the complex structure is not induced ``directly'' from $M$.

\begin{example}\label{example;notcomplex}
Let $\U(n)$ act diagonally on $p$ copies of $\C^n$.  Viewing an
element of $\C^{n\times p}$ as an $n\times p$-matrix, and identifying
$\lie u(n)^*$ with $\lie u(n)$ by means of the trace form, we can
write the moment map as $\Phi(x)=-\frac{i}{2}xx^*$.  Let $x_1$,
$x_2,\dots$, $x_n$ denote the row vectors of $x$.  The open stratum of
the imploded cross-section is the set of all $x\in\C^{n\times p}$ such
that
$$
\langle x_j,x_k\rangle=0\quad\text{for $j\ne k$,
and}\quad\norm{x_1}>\norm{x_2}>\dots>\norm{x_n}.
$$
For $p=1$ this happens to be a complex submanifold of $\C^{n\times p}$
(namely the set of vectors $(z,0,0,\dots,0)$ in $\C^n$ with $z\ne0$),
but for $p>1$ it is not.
\end{example}

Instead, the complex structure is defined indirectly, by using the
isomorphism \eqref{equation;git}.  Let $\Psi$ be the moment map for
the diagonal $K$-action on $M\times G_{\!N}$, so that
$M\impl\cong\Psi\inv(0)/K$, and let $\ca S$ be the semistable set
$$
\ca S=\{\,x\in M\times G_{\!N}\mid\text{$\overline{Gx}$ intersects
$\Psi\inv(0)$}\,\}.
$$
Results of Heinzner and Loose \cite{heinzner-loose;reduction-complex}
show that $\ca S$ is open in $M$ and that for every $x\in\ca S$ the
intersection $\overline{Gx}\cap\Psi\inv(0)$ consists of a single
$K$-orbit, so that there is a natural surjection $\ca
S\to\Psi\inv(0)/K\cong M\impl$.  Let $\ca O$ be the pushforward to
$M\impl$ of the sheaf of $G$-invariant holomorphic functions on $\ca
S$.  See e.g.\ \cite{heinzner-loose;reduction-complex} for the
definition of a K\"ahler metric on an analytic space.

\begin{theorem}\label{theorem;kaehler}
Let $M$ be a K\"ahler Hamiltonian $K$-manifold.  Then $(M\impl,\ca O)$
is an analytic space.  The strata of $M\impl$ are K\"ahler manifolds
and $M\impl$ possesses a unique K\"ahler metric which restricts to the
given K\"ahler metrics on the strata.
\end{theorem}

\begin{proof}
That the strata are K\"ahler follows from \eqref{equation;git}, which
shows that every stratum is a symplectic quotient of a K\"ahler
manifold.  The other assertions follow from Theorem (3.3) and Remark
(3.4) in \cite{heinzner-loose;reduction-complex}.
\end{proof}

\section{Quantization and implosion}\label{section;quantization}

Throughout this section let $(M,\omega,\Phi)$ be a compact Hamiltonian
$K$-manifold.  Suppose that $M$ is equivariantly prequantizable and
let $L$ be an equivariant prequantum line bundle.  By the
\emph{quantization} of $M$ we mean the equivariant index of the
Dolbeault-Dirac operator on $M$ with coefficients in $L$.  This is an
element of the representation ring of $K$ (see e.g.\
\cite{jeffrey-kirwan;localization-quantization} or
\cite{meinrenken-sjamaar;singular} for the definition), and is denoted
by $\RR(M,L)$.  In this section we compare the quantization of $M$
with that of its imploded cross-section.  A priori this does not make
sense, because $M\impl$ is not a symplectic manifold, but, following
the strategy of \cite{meinrenken-sjamaar;singular}, we shall define
the quantization of $M\impl$ to be that of a certain partial
desingularization $\tilde M\impl$.

Let $l_1$ and $l_2$ be points in $L$ with basepoints $m_1$ and $m_2$,
respectively, and define $l_1\sim l_2$ if there exists
$k\in[K_{\Phi(m_1)},K_{\Phi(m_1)}]$ such that $l_2=kl_1$.  The
\emph{imploded prequantum bundle} is the quotient $L\impl=L/{\sim}$.
There is a natural projection $L\impl\to M\impl$ and it follows from
\cite[Lemma 3.11(3)]{meinrenken-sjamaar;singular} that the fibres are
of the form $\C/\Gamma$, where $\Gamma$ is finite cyclic.  (In fact,
it is not hard to show that the restriction of $L\impl$ to each
stratum in $M\impl$ is a prequantum orbibundle.)

Let $\tau=\sigma\prin$ be the principal face of $M$.  Choose
$\lambda_0\in\tau$ and define
\begin{equation}\label{equation;desingularization}
\tilde M\impl=(M\impl\times X[\tau])\quot[\lambda_0]T\quad
\text{and}\quad\tilde L\impl=(L\impl\boxtimes\C)\quot[\lambda_0]T,
\end{equation}
where $X[\tau]$ is the symplectic toric manifold associated with the
polyhedron $-\bar\tau$ and $\C$ is the trivial line bundle on
$X[\tau]$.  In other words, $\tilde M\impl$ is the symplectic cut of
$M\impl$ with respect to the polyhedral cone $\lambda_0+\bar\tau$, and
$\tilde L\impl$ the cut bundle induced by $L\impl$.  (See
\cite{lerman;symplectic-cuts;mathematical-research} for symplectic
cutting and \cite[Section 5.2.1]{meinrenken-sjamaar;singular} for
symplectic cutting with respect to a polytope.)

Although $\tilde M\impl$ is defined as a quotient of a singular
object, observe that the fibre over $\lambda_0$ of the $T$-moment map
on $M\impl\times X[\tau]$ misses the singularities of $M\impl$.
Recall from Theorem \ref{theorem;decomposition} that the top stratum
of $M\impl$ is isomorphic as a Hamiltonian $T$-manifold to the
principal cross-section $M_\tau$.  Thus we have in fact
$$
\tilde M\impl=(M_\tau\times X[\tau])\quot[\lambda_0]T\quad
\text{and}\quad\tilde
L\impl=(L|_{M_\tau}\boxtimes\C)\quot[\lambda_0]T.
$$
This shows that for generic values of $\lambda_0$, $\tilde M\impl$ is
a Hamiltonian $T$-orbifold with moment map $\tilde\Phi\impl$, whose
image is equal to $\Phi(M)\cap(\lambda_0+\bar\tau)$.  The subset
$\tilde\Phi\impl\inv(\lambda_0+\tau)$ is a dense open submanifold,
which is isomorphic as a Hamiltonian $T$-manifold to the open
submanifold $\Phi\impl\inv(\lambda_0+\tau)$ of the top stratum
$\Phi\impl\inv(\tau)$ of $M\impl$.  Thus, as $\lambda_0$ tends to $0$,
this open set approaches the top stratum of $M\impl$.  It is in this
sense that $\tilde M\impl$ is a partial desingularization of $M\impl$,
similar to Kirwan's partial desingularization \cite{kirwan;partial} of
a symplectic quotient.  Although there is no canonical ``blow-down''
map $\tilde M\impl\to M\impl$, we shall see below in what way $\tilde
M\impl$ is a conventional desingularization of $M\impl$.

We \emph{define} the quantization of $M\impl$ to be the
$T$-equivariant index of $\tilde M\impl$ with coefficients in $\tilde
L\impl$.  In other words,
\begin{equation}\label{equation;implosion-rr}
\RR(M\impl,L\impl)=\RR(\tilde M\impl,\tilde L\impl).
\end{equation}
Now let $\Ind$ denote the holomorphic induction functor.

\begin{theorem}\label{theorem;quantization-implosion}
Let $\tau$ be the principal face of $M$ and let $\lambda_0\in\tau$ be
a sufficiently small generic element.  Then
$\RR(M,L)=\Ind_T^K\RR(M\impl,L\impl)$.  Hence quantization commutes
with implosion in the sense that $\RR(M\impl,L\impl)=\RR(M,L)^N$.
\end{theorem}

\begin{proof}
The first assertion follows immediately from \cite[Theorem
6.8]{meinrenken-sjamaar;singular} and the definition
\eqref{equation;implosion-rr}.  Taking the $N$-invariant parts of both
sides we get
\begin{equation}\label{equation;invariant-induce}
\RR(M,L)^N =(\Ind_T^KV)^N
\end{equation}
as virtual characters of $T$, where $V=\RR(\tilde M\impl,\tilde
L\impl)$.  By construction the moment polytope of $\tilde M\impl$ lies
within the fundamental chamber $\t^*_+$, so, by the quantization
commutes with reduction theorem \cite[Theorem
2.9]{meinrenken-sjamaar;singular}, the weights occurring in $V$ are
all dominant.  Let $\lambda$ be such a weight and $\C_\lambda$ the
representation of $T$ with weight $\lambda$.  By the Borel-Weil-Bott
Theorem $\Ind^K_T\C_\lambda\cong V_\lambda$, the irreducible
representation with highest weight $\lambda$.  Hence
$(\Ind^K_T\C_\lambda)^N\cong\C_\lambda$, since only the highest weight
vector $v_\lambda$ is invariant under $N$.  We conclude that
$(\Ind_T^KV)^N\cong V$.  Together with
\eqref{equation;invariant-induce} this proves the second assertion.
\end{proof}

\begin{example}
Taking $M=T^*K$ we find $\RR((T^*K)\impl,L\impl)=\RR(T^*K,L)^N$, which
by the Peter-Weyl Theorem is equal to the sum of the $V_\lambda$ over
all dominant weights $\lambda$.  Thus $(T^*K)\impl$ is a \emph{model}
for $K$ in the sense that every irreducible module occurs in its
quantization exactly once.  This application of our theorem is of
course illegal, because $T^*K$ is not compact, but the conclusion
appears correct and it would be of some interest to justify it
directly.  Cf.\ also
\cite{chuah;kaehler;transactions,chuah;kaehler;;1993}, where it is
proved that the K\"ahler quantization of the stratum of $(T^*K)\impl$
corresponding to a face $\sigma$ is the Hilbert direct sum of the
$V_\lambda$ over all $\lambda$ in $\sigma$.
\end{example}

Now assume that $M$ carries a $K$-invariant compatible complex
structure.  Then $M$ is a K\"ahler manifold and Theorem
\ref{theorem;kaehler} says that its imploded cross-section is a
K\"ahler space.  It follows from \eqref{equation;desingularization}
that the orbifold $\tilde M\impl$ is K\"ahler as well.  Following
Kirwan \cite{kirwan;partial} we call a \emph{partial
desingularization} of an analytic space $X$ any analytic orbifold
$\tilde X$ such that there exists a proper surjective bimeromorphic
map $\tilde X\to X$.

\begin{theorem}\label{theorem;kaehler-desingularization}
Let $M$ be a compact K\"ahler Hamiltonian $K$-manifold.  For
sufficiently small generic values of $\lambda_0$\upn, $\tilde M\impl$
is a partial desingularization of $M\impl$.
\end{theorem}

\begin{proof}
Let $X$ be the K\"ahler space $M\impl\times X[\tau]$, which is
equipped with a Hamiltonian $T$-action.  Denote the $T$-moment map on
$X$ by $\Psi$.  We will show that for all sufficiently small
$\mu\in\Psi(X)$ there exists a bimeromorphic map $X\quot[\mu]T\to
X\quot[0]T$.  (Properness and surjectivity are then immediate from the
compactness of $X\quot[\mu]T$ and the irreducibility of $X\quot[0]T$.)
This is well-known in the algebraic category.  Let us briefly indicate
how the argument carries over to the analytic category thanks to
results of Heinzner and Huckleberry
\cite{heinzner-huckleberry;potentials-convexity}.

Let $H=T^\C$.  The set of \emph{$\mu$-semistable points} is
$$
X\sst_\mu=\{\,x\in X\mid\text{$\overline{Hx}$ intersects
$\Psi\inv(\mu)$}\,\}.
$$
It is open and dense, if nonempty
(\cite[\S9]{heinzner-huckleberry;potentials-convexity}).  Two points
in $X\sst_\mu$ are \emph{equivalent} under the $H$-action if their
orbit closures intersect in $X\sst_\mu$.  For every $x\in X\sst_\mu$
there is a unique $y\in X\sst_\mu$ such that $Hy$ is closed in
$X\sst_\mu$ and $y$ is in the closure of $Hx$.  This implies that the
inclusion $\Psi\inv(\mu)\to X\sst_\mu$ induces a homeomorphism
$X\quot[\mu]T\to X\sst_\mu/\sim$.  (These assertions follow from the
holomorphic slice theorem,
\cite[\S2.7]{heinzner-loose;reduction-complex} or
\cite[\S0]{heinzner-huckleberry;potentials-convexity}.)  A
$\mu$-semistable point is \emph{$\mu$-stable} if $Hx$ is closed in
$X\sst_\mu$ and $H_x$ is finite.  The set of stable points is denoted
$X\st_\mu$.  It too is open and dense, if nonempty.  A point
$x\in\Psi\inv(\mu)$ is stable if and only if $T_x$ is finite.  (These
facts follow also from the holomorphic slice theorem.)  The last fact
we need is a generalization of Atiyah's result
\cite{atiyah;convexity-commuting} that for every $x\in X$ the image
$\Psi(\overline{Hx})$ is the convex hull in $\t^*$ of the $H$-fixed
points contained in $\overline{Hx}$.  Furthermore,
$\Psi(\overline{Hx})$ is equal to the full image $\Psi(X)$ for all $x$
in an open dense subset $X^\circ$.  The convexity is proved in
\cite{heinzner-huckleberry;potentials-convexity}.  For the set
$X^\circ$ we can take $X\sst_{\mu_1}\cap X\sst_{\mu_2}\cap\dots\cap
X\sst_{\mu_s}$, where $\mu_1$, $\mu_2,\dots$, $\mu_s$ are the vertices
of $\Psi(X)$.

Take $\mu\in\Psi(X)$ so small that $\Psi\inv(\mu)$ is contained in
$X\sst_0$.  Then $X\sst_\mu\subset X\sst_0$, and this inclusion
induces an analytic map $X\quot[\mu]T\cong X\sst_\mu/{\sim}\to
X\quot[0]T\cong X\sst_0/{\sim}$.  To see that this map is
bimeromorphic, observe that the stable set $X\st_0$ is nonempty, since
$0$ is a regular value of the $T$-moment map on the manifold
$M_\tau\times X[\tau]$.  Let $Y=X\st_0\cap X^\circ$.  Then the image
of $Y$ in $X\quot[0]T$ is open and dense and for $x\in Y$ we have
$\mu\in\Psi(\overline{Gx})$, i.e.\ $Y\subset X\sst_\mu$.  Thus we
obtain an analytic map $Y/{\sim}\to X\quot[\mu]T$ which inverts the
previously defined map $X\quot[\mu]T\to X\sst_0/{\sim}$ over an open
dense set.
\end{proof}

There is a more illuminating construction of this partial
desingularization for the universal imploded cross-section
$(T^*K)\impl$.  Assume that $K$ is semisimple and simply connected.
In Proposition \ref{proposition;embedding} we identified $(T^*K)\impl$
with the algebraic variety $G_{\!N}$.  We can characterize its
desingularization $(T^*K)\impl\sptilde$ in a similar manner.  Let
$\tilde G_{\!N}$ be the homogeneous vector bundle $G\times^BE^N$ over
the flag variety $G/B$ with fibre $E^N$.  Here $G=K^\C$ and $E$ is as
in \eqref{equation;bigmodule}.  The multiplication map $G\times E^N\to
E$ induces a proper morphism $p\colon\tilde G_{\!N}\to E$.

\begin{proposition}\label{proposition;isomorphism}
Suppose that $K$ is semisimple and simply connected.  Let
$\lambda_0\in\t^*$ be regular dominant and let $\omega_0$ be the
K\"ahler form on $G/B$ obtained by identifying $G/B$ with the
coadjoint $K$-orbit through $\lambda_0$.  Let $q\colon\tilde
G_{\!N}\to G/B$ be the bundle projection and put
$\tilde\omega_0=p^*\omega_E+q^*\omega_0$.
\begin{enumerate}
\item\label{item;desingularization}
$\tilde G_{\!N}$ is an equivariant desingularization of $G_{\!N}$.
\item\label{item;kahler}
$\tilde\omega_0$ is a K\"ahler form on $\tilde G_{\!N}$.  It is
integral if $\lambda_0$ is.
\item\label{item;isomorphism}
$(T^*K)\impl\sptilde$ is a smooth manifold and is isomorphic as a
Hamiltonian $K$-manifold to $\tilde G_{\!N}$.
\end{enumerate}
\end{proposition}

\begin{proof}
The image of the map $p$ is the subvariety $G_{\!N}$ of $E$ and we can
therefore regard it as a proper morphism $\tilde G_{\!N}\to G_{\!N}$.
The $G$-orbits in $\tilde G_{\!N}$ are in natural one-to-one
correspondence with the $B$-orbits in $E^N$, which are identical to
the $T^\C$-orbits in $E^N$.  Each $B$-orbit in $E^N$ passes through a
unique point of the form $v_\sigma$, so each $G$-orbit in $\tilde
G_{\!N}$ passes through a unique point of the form $[1,v_\sigma]$.
(Here points in $\tilde G_{\!N}$ are written as $[g,v]$ with $g\in G$
and $v\in E^N$.)  The stabilizer of $[1,v_\sigma]$ for the $G$-action
is $G_{[1,v_\sigma]}=B_{v_\sigma}=B\cap[P_\sigma,P_\sigma]$, where the
second equality follows from Lemma \ref{lemma;stabilizer}.  Thus the
fibre $p\inv(v_\sigma)$ is the flag variety
$[P_\sigma,P_\sigma]/(B\cap[P_\sigma,P_\sigma])$.  In particular,
$\tilde G_{\!N}$ contains a Zariski-open orbit of type $G/N$, namely
the orbit through $[1,v_\tau]$, where $\tau$ is the top face of
$\t^*_+$.  Hence $p$ is birational, which proves
\eqref{item;desingularization}.

If $\lambda_0$ is integral, then $\omega_0$ is integral on $G/B$.
Since $\omega_E$ is exact, this implies that $\tilde\omega_0$ is
integral.  Furthermore, being a sum of pullbacks of two K\"ahler
forms, $\tilde\omega_0$ is positive semidefinite.  To prove that it is
K\"ahler, it is therefore enough to show that it is nondegenerate.  We
shall do this by showing that $\tilde\omega_0$ pulls back to the
symplectic form on $(T^*K)\impl\sptilde$ under a suitable
diffeomorphism.

The principal face $\tau$ of $T^*K$ is the top face of $\t^*_+$ and
its principal cut (for the right $K$-action) is $K\times\tau$.  We
noted in Remark \ref{remark;toric} that the toric manifold associated
to the polyhedral cone $\t^*_+$ is the symplectic vector space $E^N$,
so by \eqref{equation;desingularization} the partial desingularization
of $(T^*K)\impl$ is $(K\times\tau\times E^N)\quot[\lambda_0]T$.  To
see that this space is actually a manifold rather than an orbifold,
observe that the moment map for the $T$-action on the product
$K\times\tau\times E^N$ is given by
$\Psi(k,\lambda,v)=\lambda+\Phi_{E^N}(v)$, where $\Phi_{E^N}$ is the
$T$-moment map on $E^N$.  The map $K\times E^N\to K\times\tau\times
E^N$ which sends $(k,v)$ to $(k,\lambda_0-\Phi_{E^N}(v),v)$ is a
$K\times T$-equivariant diffeomorphism onto $\Psi\inv(\lambda_0)$.  It
therefore descends to a $K$-equivariant diffeomorphism
\begin{equation}\label{equation;cut-couple}
K\times^TE^N\to(K\times\tau\times
E^N)\quot[\lambda_0]T=(T^*K)\impl\sptilde,
\end{equation}
which shows that $(T^*K)\impl\sptilde$ is smooth.  Moreover, the
inclusion map $K\to G$ induces a diffeomorphism $K\times^TE^N\to
G\times^BE^N=\tilde G_{\!N}$.  Composing with the inverse of the map
\eqref{equation;cut-couple} we obtain a diffeomorphism
$$
\tilde{\ca F}\colon(T^*K)\impl\sptilde\to\tilde G_{\!N}.
$$
To finish the proof of \eqref{item;kahler} and
\eqref{item;isomorphism} we must show that $\tilde{\ca
F}^*\tilde\omega_0$ is the symplectic form on $(T^*K)\impl\sptilde$.
Recall that the symplectic cut $(T^*K)\impl\sptilde$ contains a copy
of $K\times(\lambda_0+\tau)$ as an open dense submanifold.  The
symplectic form on this subset is the form $\omega_\tau$ of Lemma
\ref{lemma;exact} and the embedding $\ca I_0\colon
(K\times(\lambda_0+\tau))\to(T^*K)\impl\sptilde$ is given by
$$
\ca I_0(k,\lambda)=[k,\lambda,s_0(\lambda)]\in\Psi\inv(\lambda_0)/T
\subset(K\times\tau\times E^N)/T.
$$
Here $s_0\colon\lambda_0+\t^*_+\to E^N$ is any section of the map
$-\Phi_{E^N}+\lambda_0$, such as for example
$$
s_0(\lambda)=\frac1{\sqrt{\pi}}
\sum_{p=1}^r\sqrt{(\lambda-\lambda_0)(\alpha\spcheck_p)}\,v_p.
$$
Let us denote the open embedding $\tilde{\ca F}\circ\ca I_0\colon
K\times(\lambda_0+\tau)\to\tilde G_{\!N}$ by $\tilde{\ca F}_0$.  We
need to show that
\begin{equation}\label{equation;forms}
\tilde{\ca F}_0^*\tilde\omega_0=\omega_\tau.
\end{equation}
It suffices to check this identity at points of the form $(1,\lambda)$
with $\lambda\in\lambda_0+\tau$.  For $(\xi_1,\mu_1)$ and
$(\xi_2,\mu_2)$ in
$T_{(1,\lambda)}(K\times(\lambda_0+\tau))\cong\k\times\t^*$ one
readily checks that
\begin{equation}\label{equation;product-form}
(\omega_\tau)_{(1,\lambda)}((\xi_1,\mu_1)\wedge(\xi_2,\mu_2))
=\mu_1(\xi_2)-\mu_2(\xi_1)-\lambda([\xi_1,\xi_2]).
\end{equation}
On the other hand, $\tilde{\ca F}_0^*\tilde\omega_0=(p\circ\tilde{\ca
F}_0)^*\omega_E+(q\circ\tilde{\ca F}_0)^*\omega_{\lambda_0}$.  Now
$q\circ\tilde{\ca F}_0(k,\lambda)=\bar k$, where $\bar k\in K/T=G/B$
denotes the coset of $k\in K$, so
\begin{equation}\label{equation;orbit}
((q\circ\tilde{\ca F}_0)^*\omega_{\lambda_0})_{(1,\lambda)}
((\xi_1,\mu_1)\wedge(\xi_2,\mu_2))
=(\omega_{\lambda_0})_{\bar1}(\bar\xi_1\wedge\bar\xi_2)
=-\lambda_0([\xi_1,\xi_2]).
\end{equation}
A computation as in the proof of Proposition
\ref{proposition;embedding} yields
\begin{multline*}
((p\circ\tilde{\ca F}_0)^*\omega_{\lambda_0})_{(1,\lambda)}
((\xi_1,\mu_1)\wedge(\xi_2,\mu_2))\\
=\mu_1(\xi_2)-\mu_2(\xi_1)+\omega_E\bigl(\xi_{1,E}(s_0(\lambda)),
\xi_{2,E}(s_0(\lambda))\bigr),
\end{multline*}
where
\begin{align*}
\omega_E\bigl(\xi_{1,E}(s_0(\lambda)),\xi_{2,E}(s_0(\lambda))\bigr)
&=\bigl\{\Phi_E^{\xi_1},\Phi_E^{\xi_2}\bigr\}(s_0(\lambda))
=\Phi_E^{[\xi_1,\xi_2]}(s_0(\lambda))\\
&=-\frac1{2}\Im
\bigl\langle[\xi_1,\xi_2](s_0(\lambda)),s_0(\lambda)\bigr\rangle\\
&=-\frac1{2\pi}\Im2\pi i\sum_{p=1}^r
(\lambda-\lambda_0)(\alpha\spcheck_p)\,\varpi([\xi_1,\xi_2])\\
&=(\lambda_0-\lambda)([\xi_1,\xi_2]).
\end{align*}
Combining this with \eqref{equation;orbit} gives
$$
(\tilde{\ca F}_0^*\tilde\omega_0)_{(1,\lambda)}
((\xi_1,\mu_1)\wedge(\xi_2,\mu_2))
=\mu_1(\xi_2)-\mu_2(\xi_1)-\lambda([\xi_1,\xi_2]),
$$
which together with \eqref{equation;product-form} proves
\eqref{equation;forms}.
\end{proof}

Returning to the case of a K\"ahler Hamiltonian $K$-manifold $M$, let
us denote its principal face by $\tau$ and let us assume for
simplicity that $\tau$ is the top face of $\t^*_+$.  By Lemma
\ref{lemma;semisimple} we may assume that $K$ is semisimple and simply
connected.  Using Theorem \ref{theorem;universal}, Proposition
\ref{proposition;embedding} and reduction in stages we see that
\begin{equation}\label{equation;universal-desingularization}
\begin{split}
\tilde M\impl=(M\impl\times X[\tau])\quot[\lambda_0]T
&\cong\bigl((M\times G_{\!N})\quot[0]K\times
X[\tau]\bigr)\bigquot[\lambda_0]T\\
&\cong\bigl(M\times(G_{\!N}\times
X[\tau])\quot[\lambda_0]T\bigr)\bigquot[0]K\\
&\cong\bigl(M\times\tilde G_{\!N}\bigr)\bigquot[0]K,
\end{split}
\end{equation}
because $\tau$ is also the principal face of $G_{\!N}$.  Thus the
desingularization of $G_{\!N}$ plays a universal role analogous to
that of $G_{\!N}$ itself.  By Proposition
\ref{proposition;isomorphism}\eqref{item;desingularization} the
analytic map
$$
p_M=\id_M\times p\colon M\times\tilde G_{\!N}\to M\times G_{\!N},
$$
is a $G$-equivariant desingularization.  The following is now clear
from Theorem \ref{theorem;kaehler-desingularization}.

\begin{corollary}
Let $M$ be a compact K\"ahler Hamiltonian $K$-manifold.  Assume that
the principal face of $M$ is the top face of $\t^*_+$.  Then for small
generic $\lambda_0$ the map $p_M$ induces a map $\tilde M\impl\to
M\impl$ which is identical to the partial desingularization of Theorem
{\rm\ref{theorem;kaehler-desingularization}}.
\end{corollary}

We mention without proof that a similar result holds if the principal
face $\tau$ is not the top face.  In this case the toric manifold
$X[\tau]$ is not $E^N$, but the smaller symplectic vector space
$E^{[P_\tau,P_\tau]}$.  In Proposition \ref{proposition;isomorphism}
and in \eqref{equation;universal-desingularization} one needs to
replace $G_{\!N}$ by the closure of the stratum $G/[P_\tau,P_\tau]$
and $\tilde G_{\!N}$ by the homogeneous bundle
$G\times^{P_\tau}E^{[P_\tau,P_\tau]}$ over the partial flag variety
$G/P_\tau$.



\section{Notation index}\label{section;notation}

\begin{tabbing}
\indent \= $M_\mu$; $M_0=M\quot K$ \= $100$ \= \kill 
\> $K$; $T$ \> compact connected Lie group; maximal torus,
\S\ref{section;construction}\\
\> $Z$ \> identity component of centre, \S\ref{section;construction}\\
\> $R$; $R_+$; $S$ \> root system; positive roots; simple roots,
\S\ref{section;stratification}\\
\> $G$ \> complexification of $K$, \S\ref{section;kaehler}\\
\> $N$; $B$ \> maximal unipotent subgroup; Borel subgroup,
\S\ref{section;kaehler}\\
\> $\t^*_+$; $\Sigma$ \> closed Weyl chamber in $\lie t^*$; set of
open faces of $\t^*_+$, \S\ref{section;construction}\\
\> $\Lambda^*$; $\Lambda^*_+$ \> weight lattice in $\lie t^*$; set of
dominant weights, \S\ref{section;kaehler}\\
\> $K_\sigma$; $R_\sigma$ \> centralizer of face $\sigma\in\Sigma$;
root system of $K_\sigma$, \S\ref{section;construction}\\
\> $P_\sigma$; $U_\sigma$ \> parabolic associated to $\sigma$;
unipotent radical of $P_\sigma$, \S\ref{section;kaehler}\\
\> $G_\sigma$; $N_\sigma$ \> Levi factor of $P_\sigma$; maximal
unipotent subgroup of $G_\sigma$, \S\ref{section;kaehler}\\
\> $\eu S_\sigma$\> standard slice at $\sigma$ for coadjoint action,
\S\ref{section;construction}\\
\> $X_N$ \> variety with coordinate ring $\C[X]^N$,
\S\ref{section;kaehler}\\
\> $(M,\omega)$; $\Phi$ \> Hamiltonian $K$-manifold; moment map,
\S\ref{section;construction}\\
\> $M\quot[\lambda]K$; $M\quot K$ \> symplectic quotient
$\Phi\inv(\lambda)/K_\lambda$; same with $\lambda=0$,
\S\ref{section;construction}\\
\> $M_\sigma$; $\sigma\prin$ \> cross-section $\Phi\inv(\eu
S_\sigma)$; principal face of $M$, \S\ref{section;construction}\\
\> $\sim$; $M\impl$ \> equivalence relation on $\cross$; imploded
cross-section, \S\ref{section;construction}\\
\> $\pi$; $\Phi\impl$ \> quotient map $\cross\to M\impl$; imploded
moment map, \S\ref{section;construction}\\
\> $\ca L$; $\ca R$ \> left; resp.\ right action of $K$ on itself or
on $T^*K$, \S\ref{section;universal}\\
\> $\pi_{\ca R}$ \> quotient map $\Phi_{\ca
R}\inv(\t^*_+)\to(T^*K)\impl$, \S\ref{section;universal}\\
\> $\RR(M,L)$ \> equivariant index of $M$ with coefficients in $L$,
\S\ref{section;quantization}\\
\> $C^\circ(X)$ \> infinite cone
$\bigl(X\times[\,0,\infty)\bigr)\big/\bigl(X\times\{0\}\bigr)$,
\S\ref{section;universal}\\
\> $W^\omega$ \> skew complement, \S\ref{section;stratification}
\end{tabbing}


\providecommand{\bysame}{\leavevmode\hbox to3em{\hrulefill}\thinspace}
\providecommand{\MR}{\relax\ifhmode\unskip\space\fi MR }
\providecommand{\MRhref}[2]{%
  \href{http://www.ams.org/mathscinet-getitem?mr=#1}{#2}
}
\providecommand{\href}[2]{#2}



\begin{thebibliography}{10}

\bibitem{arnold;mathematical-methods}
V.~Arnold, \emph{Mathematical methods of classical mechanics}, corrected
  reprint of second {E}nglish (1989) ed., Graduate Texts in Mathematics,
  vol.~60, Springer-Verlag, New York, 199?, translated from the Russian by K.
  Vogtmann and A. Weinstein. \MR{96c:70001}

\bibitem{atiyah;convexity-commuting}
M.~Atiyah, \emph{Convexity and commuting {H}amiltonians}, Bull. London Math.
  Soc. \textbf{14} (1982), no.~1, 1--15. \MR{83e:53037}

\bibitem{bernstein-gelfand-gelfand;differential-operators-base}
I.~Bernstein, I.~Gelfand, and S.~Gelfand, \emph{Differential operators on the
  base affine space and a study of ${\mathfrak g}$-modules}, Lie Groups and
  their Representations (Budapest, 1971) (I.~Gelfand, ed.), Bolyai J\'anos
  Math. Soc., Halsted, New York, 1975, pp.~21--64. \MR{58 \#28285}

\bibitem{chuah;kaehler;transactions}
M.-K. Chuah, \emph{Kaehler structures on {$K_{\mathbf C}/(P,P)$}}, Trans. Amer.
  Math. Soc. \textbf{349} (1997), no.~8, 3373--3390. \MR{97k:22016}

\bibitem{chuah;kaehler;;1993}
M.-K. Chuah and V.~Guillemin, \emph{Kaehler structures on {${K}_{\mathbf
  C}/N$}}, The Penrose transform and analytic cohomology in representation
  theory (South Hadley, MA, 1992) (M.~Eastwood et~al., eds.), Contemporary
  Mathematics, vol. 154, Amer. Math. Soc., Providence, RI, 1993, pp.~181--195.
  \MR{94k:22028}

\bibitem{guillemin-lerman-sternberg;symplectic-fibrations}
V.~Guillemin, E.~Lerman, and S.~Sternberg, \emph{Symplectic fibrations and
  multiplicity diagrams}, Cambridge University Press, Cambridge, 1996.
  \MR{98d:58074}

\bibitem{heinzner-huckleberry;potentials-convexity}
P.~Heinzner and A.~Huckleberry, \emph{K\"ahlerian potentials and convexity
  properties of the moment map}, Invent. Math. \textbf{126} (1996), no.~1,
  65--84. \MR{98e:58075}

\bibitem{heinzner-loose;reduction-complex}
P.~Heinzner and F.~Loose, \emph{Reduction of complex {H}amiltonian
  ${G}$-spaces}, Geom. Funct. Anal. \textbf{4} (1994), no.~3, 288--297.
  \MR{95j:58050}

\bibitem{hurtubise-jeffrey;representations-weighted}
J.~Hurtubise and L.~Jeffrey, \emph{Representations with weighted frames and
  framed parabolic bundles}, Canad. J. Math. \textbf{52} (2000), 1235--1268.

\bibitem{jeffrey-kirwan;localization-quantization}
L.~Jeffrey and F.~Kirwan, \emph{Localization and the quantization conjecture},
  Topology \textbf{36} (1997), no.~3, 647--693. \MR{98e:58088}

\bibitem{kempf-ness;length}
G.~Kempf and L.~Ness, \emph{The length of vectors in representation spaces},
  Algebraic geometry (Copenhagen, 1978) (K.~L{\o}nsted, ed.), Lecture Notes in
  Mathematics, vol. 732, Springer-Verlag, Berlin, 1979, pp.~233--243.
  \MR{81i:14032}

\bibitem{kirwan;partial}
F.~Kirwan, \emph{Partial desingularisations of quotients of nonsingular
  varieties and their {B}etti numbers}, Ann. of Math. (2) \textbf{122} (1985),
  no.~1, 41--85. \MR{87a:14010}

\bibitem{kraft;geometrische}
H.~Kraft, \emph{Geometrische {M}ethoden in der {I}nvariantentheorie}, second
  revised ed., Aspekte der {M}athematik, vol.~D1, Vieweg, Braunschweig, 1985.
  \MR{86j:14006}

\bibitem{lerman;symplectic-cuts;mathematical-research}
E.~Lerman, \emph{Symplectic cuts}, Math. Res. Lett. \textbf{2} (1995), no.~3,
  247--258. \MR{96f:58062}

\bibitem{lerman;nonabelian-convexity}
E.~Lerman, E.~Meinrenken, S.~Tolman, and C.~Woodward, \emph{Nonabelian
  convexity by symplectic cuts}, Topology \textbf{37} (1998), no.~2, 245--259.
  \MR{99a:58069}

\bibitem{mcduff-salamon;introduction}
D.~McDuff and D.~Salamon, \emph{Introduction to symplectic topology}, Oxford
  Mathematical Monographs, Oxford University Press, New York, 1995.
  \MR{97b:58062}

\bibitem{meinrenken-sjamaar;singular}
E.~Meinrenken and R.~Sjamaar, \emph{Singular reduction and quantization},
  Topology \textbf{38} (1999), no.~4, 699--762. \MR{2000f:53114}

\bibitem{pauer;glatte}
F.~Pauer, \emph{Glatte {E}inbettungen von {$G/U$}}, Math. Ann. \textbf{262}
  (1983), no.~3, 421--429. \MR{84i:14032}

\bibitem{prill;local}
D.~Prill, \emph{Local classification of quotients of complex manifolds by
  discontinuous groups}, Duke Math. J. \textbf{34} (1967), 375--386. \MR{35
  \#1829}

\bibitem{schwarz;topology-algebraic}
G.~Schwarz, \emph{The topology of algebraic quotients}, Topological methods in
  algebraic transformation groups (Rutgers Univ., 1988) (H.~Kraft et~al.,
  eds.), Progress in Mathematics, vol.~80, Birkh\"auser, Boston, 1989,
  pp.~135--151. \MR{90m:14043}

\bibitem{sjamaar-lerman;stratified}
R.~Sjamaar and E.~Lerman, \emph{Stratified symplectic spaces and reduction},
  Ann. of Math. (2) \textbf{134} (1991), no.~2, 375--422. \MR{92g:58036}

\bibitem{vinberg-popov;class-quasihomogeneous}
{\`E}.~B. Vinberg and V.~L. Popov, \emph{A certain class of quasihomogeneous
  affine varieties}, Izv. Akad. Nauk SSSR Ser. Mat. \textbf{36} (1972),
  749--764. \MR{47 \#1815}

\end{thebibliography}
\end{document}